\documentclass[12pt]{article}%
\usepackage{graphicx}
\usepackage[intlimits]{amsmath}
\usepackage{latexsym}
\usepackage{amsfonts}
\usepackage{amssymb}%
\makeatletter
\newfont{\footsc}{cmcsc10 at 8truept}
\newfont{\footbf}{cmbx10 at 8truept}
\newfont{\footrm}{cmr10 at 10truept}
\newtheorem{theorem}{Theorem}[section]

\newtheorem{claim}[theorem]{Claim}

\newtheorem{conjecture}[theorem]{Conjecture}
\newtheorem{corollary}[theorem]{Corollary}

\newtheorem{fact}[theorem]{Fact}
\newtheorem{lemma}[theorem]{Lemma}

\newtheorem{question}[theorem]{Question}
\newtheorem{proposition}[theorem]{Proposition}

\newenvironment{proof}[1][Proof]{\noindent{\textbf {#1}  }}  {\hfill$\Box$\bigskip}

\def\blfootnote{\xdef\@thefnmark{}\@footnotetext}
\begin{document}

\title{Ramsey Goodness and Beyond}
\author{V. Nikiforov, C. C. Rousseau\\{\small Department of Mathematical Sciences, University of Memphis, Memphis,
TN 38152}\\{\small e-mail:} {\small vnikifrv@memphis.edu, ccrousse@memphis.edu}}
\maketitle

\begin{abstract}
In a seminal paper from 1983, Burr and Erd\H{o}s started the systematic study
of \ Ramsey numbers of cliques vs. large sparse graphs, raising a number of
problems. In this paper we develop a new approach to such Ramsey problems
using a mix of the Szemer\'{e}di regularity lemma, embedding of sparse graphs,
Tur\'{a}n type stability, and other structural results. We give exact Ramsey
numbers for various classes of graphs, solving all but one of the
Burr-Erd\H{o}s problems.\medskip

\xdef\@thefnmark{}\@footnotetext{\textit{2000 Mathematics Subject Classification:
}Primary 05C55; Secondary 05C35.}

\xdef\@thefnmark{}\@footnotetext{\textit{Key words and phrases: }Ramsey numbers of sparse graphs;
Ramsey goodness; degenerate graphs; wheels; joints.}\newpage

\end{abstract}
\tableofcontents

\section{Introduction}

Our notation is standard (e.g., see \cite{Bol98}). In particular, $G\left(
n\right)  $ stands for a graph of order $n;$ we write $\left\vert G\right\vert
$ for the order of a graph $G$ and $k_{r}\left(  G\right)  $ for the number of
its $r$-cliques. The join of the graphs $G$ and $H$ is denoted by $G+H.$

Given a graph $G,$ a\emph{ }$2$\emph{-coloring} of $E\left(  G\right)  $ is a
partition $E\left(  G\right)  =E\left(  R\right)  \cup E\left(  B\right)  ,$
where $R$ and $B$ are graphs with $V\left(  R\right)  =V\left(  B\right)
=V\left(  G\right)  .$ The \emph{Ramsey number} $r(H_{1},H_{2})$ is the least
number $n$ such that for every $2$-coloring $E\left(  K_{n}\right)  =E\left(
R\right)  \cup E\left(  B\right)  ,$ either $H_{1}\subset R$ or $H_{2}\subset
B$.

The aim of this paper is to develop a new approach to Ramsey numbers of
cliques vs. large sparse graphs. We prove a generic Ramsey result about
certain classes of graphs, thus producing an unlimited source of specific
exact Ramsey numbers. This enables us to answer a number of open questions and
extend a substantial amount of earlier research. Moreover, some of the
auxiliary results used in our proofs may be regarded as general tools for
wider classes of Ramsey problems.

Let us recall the notion of goodness in Ramsey theory, introduced by Burr
\cite{Bur81}: a connected graph $H$ is $p$-good if the Ramsey number
$r(K_{p},H)$ is given by
\[
r\left(  K_{p},H\right)  =\left(  p-1\right)  \left(  \left\vert H\right\vert
-1\right)  +1.
\]

The systematic study of good Ramsey results was initiated by Burr and
Erd\H{o}s in \cite{BuEr83}; for surveys of subsequent progress the reader is
referred to \cite{Bur87} and \cite{FRS89}.

First we outline some of the problems raised in \cite{BuEr83}.

\subsection{\label{Spro}Solved and unsolved problems about $p$-good graphs}

In \cite{BuEr83} Burr and Erd\H{o}s, probing the limits of $p$-goodness, gave
some general constructions of $p$-good graphs and raised a number of
questions, most of which are still open. To state the most important problem
raised in \cite{BuEr83} and reiterated in \cite{Bur87} and \cite{EFRS85}, we
recall that a graph is called $q$\emph{-degenerate} if each of its subgraphs
contains a vertex of degree at most $q.$

\begin{conjecture}
\label{con1}For fixed $q\geq1,$ $p\geq3$, all sufficiently large
$q$-degenerate graphs are $p$-good.
\end{conjecture}

A weaker version of this conjecture was stated earlier by Burr in \cite{Bur81}.

\begin{conjecture}
\label{con2}For fixed $q\geq1,$ $p\geq3$, all sufficiently large graphs of
maximum degree at most $q$ are $p$-good.
\end{conjecture}

Brandt \cite{Bra96} showed that for $p=3$ and $q\geq168,$ every $q$-regular
graph of sufficiently large order and with sufficiently large expansion factor
is a counterexample to Conjecture \ref{con2}. Using a different approach, in
Section \ref{Discon} we show that, for $p=3,$ almost all $100$-regular graphs
are counterexamples to Conjecture \ref{con2}, and thus to Conjecture
\ref{con1}.

We shall answer in the affirmative all but one of the remaining questions
raised in \cite{BuEr83}.

Write $C_{n}$ for the cycle of order $n.$ Burr and Erd\H{o}s \cite{BuEr83}
showed that the wheel $K_{1}+C_{n}$ is $3$-good for $n\geq5.$ This result
motivated the following three questions (\cite{BuEr83}, p. 50.)

\begin{question}
\label{pr1}Is the wheel $K_{1}+C_{n}$ $p$-good for fixed $p>3$ and $n$ large?
\end{question}

Recall that the $k$th power of a graph $G$ is a graph $G^{k}$ with $V\left(
G^{k}\right)  =V\left(  G\right)  $ and $uv\in E\left(  G^{k}\right)  $ if $u$
and $v$ can be joined in $G$ by a path of length at most $k.$

\begin{question}
\label{pr2} Is $K_{1}+C_{n}^{k}$ a $p$-good graph for fixed $k\geq2,$ $p\geq3$
and $n$ large?
\end{question}

\begin{question}
\label{pr3}Fix $p\geq3,$ $l\geq1,$ $k\geq1,$ and a connected graph $G$. Is it
true that, for every large enough graph $G_{1}$ homeomorphic to $G,$ the graph
$K_{l}+G_{1}^{k}$ is $p$-good?
\end{question}

Burr and Erd\H{o}s estimated that finding an answer to Question \ref{pr3}
would be very difficult. In this paper we answer Question \ref{pr3} in the
affirmative, implying an affirmative answer to Questions \ref{pr1} and
\ref{pr2} as well.

Clearly, the clique number of $p$-good graphs must grow rather slowly with
their order. Therefore, the following question comes naturally (\cite{BuEr83},
p. 41.)

\begin{question}
\label{pr4}Subdivide each edge of $K_{n}$ by one vertex. Is the resulting
graph $p$-good for $p$ fixed and $n$ large?
\end{question}

Burr and Erd\H{o}s also asked a question about tree-like constructions of
fixed families of graphs, which they called \textquotedblleft graphs with
bridges\textquotedblright\ (\cite{BuEr83}, p. 44). We restate their question
in a much stronger form.

Given a graph $G$ of order $n$ and a vector of positive integers
$\mathbf{k}=\left(  k_{1},\ldots,k_{n}\right)  $, write $G^{\mathbf{k}}$ for
the graph obtained from $G$ by replacing each vertex $i\in\left[  n\right]  $
with a clique of order $k_{i}$ and every edge $ij\in E\left(  G\right)  $ with
a complete bipartite graph $K_{k_{i},k_{j}}.$

\begin{question}
\label{pr6}Suppose $K\geq1$, $p\geq3$, $T_{n}$ is a tree of order $n$, and
$\mathbf{k}=\left(  k_{1},\ldots,k_{n}\right)  $ is a vector of integers with
$0<k_{i}\leq K$ for all $i\in\left[  n\right]  .$ Is $T_{n}^{\mathbf{k}}$
$p$-good for $n$ large?
\end{question}

We shall answer Questions \ref{pr4} and \ref{pr6} in the affirmative. However,
the following particular question raised in \cite{BuEr83} is beyond the scope
of this paper.

\begin{question}
\label{prcub}Is the $n$-cube $3$-good for $n$ large?
\end{question}

\subsection{Some highlights on Ramsey goodness}

We list below several important results on Ramsey goodness.

Define a $q$\emph{-book} of size $n$ to be the graph $B_{q}\left(  n\right)
=K_{q}+nK_{1},$ i.e., $B_{q}\left(  n\right)  $ consists of $n$ distinct
$\left(  q+1\right)  $-cliques sharing a $q$-clique.

\begin{fact}
[\cite{NiRo04}]\label{f1}For fixed $q\geq2,$ $p\geq3,$ and large $n,$
\[
r\left(  K_{p},B_{q}\left(  n\right)  \right)  =\left(  p-1\right)  \left(
n+q-1\right)  +1.
\]

\end{fact}

In the following results $K_{p}$ is replaced by a supergraph $H\supset K_{p}$
such that $r\left(  H,G\right)  =r\left(  K_{p},G\right)  $ for certain
$p$-good graphs $G$.

\begin{fact}
[\cite{FRS78}, \cite{FRS91}]\label{f5}For fixed $m\geq1$ and large $n,$
\[
r\left(  B_{2}\left(  m\right)  ,C_{n}\right)  =2\left(  n-1\right)  +1.
\]

\end{fact}

\begin{fact}
[\cite{RoSh78b}, \cite{FSR80}]\label{f6}For fixed $p\geq2,$ $m\geq1,$ and any
tree $T_{n}$ of large order $n,$
\[
r\left(  B_{p}\left(  m\right)  ,T_{n}\right)  =p\left(  n-1\right)  +1.
\]

\end{fact}

Write $K_{p}\left(  t_{1},\ldots,t_{p}\right)  $ for the complete $p$-partite
graph with part sizes $t_{1},\ldots,t_{p}$ and set $K_{p}\left(  t\right)
=K_{p}\left(  t,\ldots,t\right)  .$

\begin{fact}
[\cite{BFRS83}, \cite{EFRS88a}]\label{f7}For fixed $m\geq1,$ $k\geq1,$
$n_{1},\ldots,n_{k},$ and $n$ large,
\[
r\left(  B_{2}\left(  m\right)  ,K_{k+1}\left(  n_{1},\ldots,n_{k},n\right)
\right)  =2\left(  n_{1}+\ldots+n_{k}+n\right)  -1.
\]

\end{fact}

\begin{fact}
[\cite{B-al87}, \cite{EFRS88a}]\label{f2}For fixed $p\geq2,$ $t\geq1,$ and any
tree $T_{n}$ of large order $n,$
\[
r\left(  K_{p+1}\left(  1,1,t,\ldots,t\right)  ,T_{n}\right)  =p\left(
n-1\right)  +1.
\]

\end{fact}

\begin{fact}
[\cite{RoSh78a}, \cite{FRS82}, \cite{NiRo05a}, \cite{NiRo05b}]\label{f3}There
exists $c>0$ independent of $n$ such that if $n$ is large and $m\leq cn$,
then
\[
r\left(  B_{2}\left(  m\right)  ,B_{2}\left(  n\right)  \right)  =2n+3.
\]

\end{fact}

The following result answers in the affirmative a special case of Question
\ref{pr6}.

\begin{fact}
[\cite{LiRo96}]\label{f4}For fixed $p\geq3$ and graph $H$, the graph
$K_{1}+nH$ is $p$-good for $n$ large.
\end{fact}

\section{Main results}

We first outline the approach to Ramsey numbers adopted in this paper.

For every $p$ and $n,$ we describe two families of graphs $\mathcal{R}\left(
n\right)  $ and $\mathcal{B}\left(  n\right)  $ such that, if $n$ is large,
then for every $2$-coloring $E\left(  K_{p\left(  n-1\right)  +1}\right)
=E\left(  R\right)  \cup E\left(  B\right)  ,$ either $H\subset R$ for some
$H\in\mathcal{R}\left(  n\right)  $ or $G\subset B$ for all $G\in
\mathcal{B}\left(  n\right)  .$

To describe $\mathcal{R}\left(  n\right)  ,$ we define joints: call the union
of $t$ distinct $p$-cliques sharing an edge a\emph{ }$p$\emph{-joint }of size
$t;$ denote the maximum size of a $p$-joint in a graph $G$ by $js_{p}\left(
G\right)  .$ The family $\mathcal{R}\left(  n\right)  $ consists of all
$\left(  p+1\right)  $-joints of size at least $cn^{p-1}$ for some appropriate
$c>0.$

To describe $\mathcal{B}\left(  n\right)  ,$ we first define splittable
graphs: given real numbers $\gamma,\eta>0,$ we say that a graph $G=G\left(
n\right)  $ is $\left(  \gamma,\eta\right)  $\emph{-splittable }if there
exists a set $S\subset V\left(  G\right)  $ with $\left\vert S\right\vert
<n^{1-\gamma}$ such that the order of any component of $G-S$ is at most $\eta
n.$ The family $\mathcal{B}\left(  n\right)  $ consists of all $q$-degenerate
$\left(  \gamma,\eta\right)  $-splittable graphs, where $q$ and $\gamma$ are
fixed and $\eta>0$ is appropriately chosen.

\subsection{The main theorem}

Here is our main theorem.

\begin{theorem}
\label{mth}For all $p\geq3,$ $q\geq1,$ $0<\gamma<1,$ there exist
$c>0,\ \eta>0$ such that if $E\left(  K_{p\left(  n-1\right)  +1}\right)
=E\left(  R\right)  \cup E\left(  B\right)  $ is a $2$-coloring, then for $n$
large, one of the following conditions holds:

(i) $R$ contains a $\left(  p+1\right)  $-joint of size $cn^{p-1};$

(ii) $B$ contains every $q$-degenerate $\left(  \gamma,\eta\right)
$-splittable graph $G$ of order $n$.
\end{theorem}

Note that Theorem \ref{mth} gives exact Ramsey numbers for graphs of varying
structure, implying, in particular, positive answers to the questions raised
in Section \ref{Spro}.

\subsection{Variations of the $\mathcal{R}\left(  n\right)  $ family}

The condition $js_{p+1}\left(  R\right)  >cn^{p-1}$ implies the existence of
various $\left(  p+1\right)  $-partite graphs in $R$. On the one hand, $R$
contains dense supergraphs of $K_{p+1}$ as shown in the following theorem,
proved in \ref{pmth1}.

\begin{theorem}
\label{mth1} For all $p\geq3,$ $q\geq1,$ $0<\gamma<1,$ there exist
$c>0,\ \eta>0$ such that if $E\left(  K_{p\left(  n-1\right)  +1}\right)
=E\left(  R\right)  \cup E\left(  B\right)  $ is a $2$-coloring, then for $n$
large, one of the following conditions holds:

(i) $R$ contains $K_{p+1}\left(  1,1,t,\ldots,t\right)  $ for $t=\left\lceil
c\log n\right\rceil ;$

(ii) $B$ contains every $q$-degenerate, $\left(  \gamma,\eta\right)
$-splittable graph $G$ of order $n$.
\end{theorem}

Observe that this theorem considerably changes the usual setup for goodness
results: now in the graph $R$ we find dense supergraphs of $K_{p}$ whose order
grows with $n.$ On the other hand, if we give up density, we find in $R$
sparse $p$-partite graphs whose order is linear in $n$. More precisely, we
have the following theorem, proved in \ref{pmth2}.

\begin{theorem}
\label{mth2} For all $p\geq2,$ $q\geq1,$ $d\geq2,$ $0<\gamma<1,$ there exist
$\alpha>0,$ $c>0,\ \eta>0$ such that if $E\left(  K_{p\left(  n-1\right)
+1}\right)  =E\left(  R\right)  \cup E\left(  B\right)  $ is a $2$-coloring,
then, for $n$ large, one of the following conditions holds:

(i) $R$ contains $K_{s}+H$ for every $\left(  p+1-s\right)  $-partite graph
$H$ with $\left\vert H\right\vert =\left\lfloor \alpha n\right\rfloor $ and
$\Delta\left(  H\right)  \leq d;$

(ii) $B$ contains every $q$-degenerate, $\left(  \gamma,\eta\right)
$-splittable graph $G$ of order $n$.
\end{theorem}

\subsection{Variations of the $\mathcal{B}\left(  n\right)  $ family}

Call a family of graphs $\mathcal{F}$ $\gamma$-\emph{crumbling}, if for any
$\eta>0$, there exists $n_{0}\left(  \eta\right)  $ such that all graphs
$G\in\mathcal{F}$ with $\left\vert G\right\vert >n_{0}\left(  \eta\right)  $
are $\left(  \gamma,\eta\right)  $-splittable. We will say that a family
$\mathcal{F}$ is degenerate and crumbling if $\mathcal{F}$ is $q$-degenerate
and $\gamma$-crumbling for some specific $q$ and $\gamma$.

Restricting Theorem \ref{mth} to degenerate crumbling families, we obtain the
following theorem.

\begin{theorem}
\label{mth3}For all $p\geq2,$ $q\geq1,$ $0<\gamma<1$ there exists $c>0$ such
that if $\mathcal{F}$ is a $q$-degenerate $\gamma$-crumbling family and
$E\left(  K_{p\left(  n-1\right)  +1}\right)  =E\left(  R\right)  \cup
E\left(  B\right)  $ is a $2$-coloring, then for $n$ large one of the
following conditions holds:

(i) $R$ contains a $\left(  p+1\right)  $-joint of size $cn^{p-1};$

(ii) $B$ contains every $G\in\mathcal{F}$ of order $n.$
\end{theorem}

Since $K_{p}$ is a subgraph of any $p$-joint, it follows that all sufficiently
large members of a degenerate crumbling family are $p$-good. This simple
observation is a clue to the answers of all questions of Section \ref{Spro}.

Subdivide each edge of $K_{n}$ by a single vertex, write $\widehat{K_{n}}$ for
the resulting graph, and note that $\widehat{K_{n}}$ is $2$-degenerate. If we
remove the vertices of the original $K_{n},$ the remaining graph consists of
$\binom{n}{2}$ isolated vertices. Since $n<\left(  n\left(  n+1\right)
\right)  ^{1/2}$ it follows that $\widehat{K_{n}}$ is $\left(  1/2,\eta
\right)  $-splittable for $\eta=1/\binom{n+1}{2}$. Thus, the family of all
$\widehat{K_{n}}$'s is $2$-degenerate and crumbling; hence, $\widehat{K_{n}}$
is $p$-good for $n$ large, answering Question \ref{pr4}.

The propositions stated below are proved in Section \ref{Degs} unless their
proof is omitted.

The answer to Question \ref{pr3} is affirmative in view of the following three propositions.

\begin{proposition}
\label{p1}The family of all graphs homeomorphic to a fixed connected graph $G$
is degenerate and crumbling.
\end{proposition}

\begin{proposition}
\label{p4}If $\mathcal{F}$ is a crumbling family of bounded maximum degree,
then, for fixed $k\geq1,$ the family $\mathcal{F}^{k}=\left\{  G^{k}%
:G\in\mathcal{F}\right\}  $ is degenerate and crumbling.
\end{proposition}

The following proposition is obvious, so we omit its proof.

\begin{proposition}
\label{p2}Let $l\geq l$ be a fixed integer. If $\mathcal{F}$ is a degenerate
crumbling family, then the family of connected graphs $\mathcal{F}^{\ast
}=\left\{  K_{l}+G:G\in\mathcal{F}\right\}  $ is degenerate and
crumbling.\hfill{$\Box$}
\end{proposition}

Note also that, in view of Proposition \ref{p2}, Theorem \ref{mth3}
generalizes Fact \ref{f4}.

Trees provide various examples of degenerate crumbling families.

\begin{proposition}
\label{p3}Every infinite family of trees is degenerate and crumbling.
\end{proposition}

In particular, Proposition \ref{p3} and Theorem \ref{mth3} extend Fact
\ref{f2}. Likewise, Theorem \ref{mth3} and the following simple observation,
whose proof we omit, extend Fact \ref{f1}.

\begin{proposition}
\label{p7}Every infinite family of $q$-books is degenerate and crumbling.
\hfill{$\Box$}
\end{proposition}

Some operations on graphs fit well with degenerate and crumbling families, as
shown in Proposition \ref{p4} and the following two propositions.

\begin{proposition}
\label{p5}Let $\mathcal{F}_{1}$ and $\mathcal{F}_{2}$ be degenerate crumbling
families. Then the family
\[
\mathcal{F}_{1}\times\mathcal{F}_{2}=\left\{  G_{1}\times G_{2}:G_{1}%
\in\mathcal{F}_{1},G_{2}\in\mathcal{F}_{2}\right\}
\]
is degenerate and crumbling.
\end{proposition}

\begin{proposition}
\label{p6}Let $\mathcal{F}$ be a degenerate crumbling family, and $\left\{
\mathbf{k}_{n}=\left(  k_{1},\ldots,k_{n}\right)  \right\}  _{n=1}^{\infty}$
be a sequence of integer vectors with $0<k_{i}\leq K,$ for $i\in\left[
n\right]  .$ Then the family $\mathcal{F}^{\ast}=\left\{  G^{\mathbf{k}_{n}%
}:G\in\mathcal{F},\left\vert G\right\vert =n\right\}  $ is degenerate and crumbling.
\end{proposition}

Note that Proposition \ref{p6}, together with Theorem \ref{mth3} answers
Question \ref{pr6} in the affirmative.

As an additional application consider the following example: write
$Grid_{n}^{k}$ for the product of $k$ copies of the path $P_{n},$ i.e.,
$V\left(  Grid_{n}^{k}\right)  =\left[  n\right]  ^{k}$ and two vertices
$\left(  u_{1},\ldots,u_{k}\right)  ,\left(  v_{1},\ldots,v_{k}\right)
\in\left[  n\right]  ^{k}$ are joined if $\sum_{i=1}^{k}\left\vert u_{i}%
-v_{i}\right\vert =1$. Propositions \ref{p3}, \ref{p5}, and Theorem \ref{mth3}
imply that $Grid_{n}^{k}$ is $p$-good for $k$ fixed and $n$ large; it seems
that this natural problem hasn't been raised earlier.

A particular instance of Theorem \ref{mth3} is the following extension of
Facts \ref{f5}, \ref{f6}, \ref{f7}, and \ref{f3}.

\begin{theorem}
\label{th2}For all $p\geq2,$ $q\geq1,$ $\gamma>0$ there exist $c>0$ such that
for every $q$-degenerate $\gamma$-crumbling family $\mathcal{F}$ of connected
graphs, then%
\[
r\left(  B_{p}\left(  \left\lceil cn\right\rceil \right)  ,G\right)  =p\left(
n-1\right)  +1
\]
for every $G\in\mathcal{F}$ of sufficiently large order $n.$
\end{theorem}

Indeed, it suffice to note that if $js_{p+1}\left(  R\right)  >cn^{p-1},$ then
$B_{p}\left(  \left\lceil cn\right\rceil \right)  \subset R.$

Restricting Theorem \ref{mth2} to crumbling degenerate families, we can
substantially generalize Theorem \ref{th2} and replace the graph $B_{p}\left(
\left\lceil cn\right\rceil \right)  $ with other graphs, e.g., $K_{p-1}%
+C_{\left\lceil cn\right\rceil },$ where $\left\lceil cn\right\rceil $ is even.

\subsection{Remarks on the proof methods}

The proof of Theorem \ref{mth} is based on several major results. The key
element is a compound of the Szemer\'{e}di regularity lemma and a structural
theorem in \cite{Nik06}, stating that, for sufficiently small $c>0,$ the
vertices of any graph $G$ with $k_{p}\left(  G\right)  <cn^{p}$ can be
partitioned into bounded number of very sparse sets. Other ingredients are a
stability result about large $p$-joints, proved in \cite{BoNi04}, and a
probabilistic lemma used in different forms by other researchers. Finally we
construct several rather involved embedding algorithms for degenerate
splittable graphs.

\section{Proofs}

We start with some additional notation. Set $\left[  n\right]  =\left\{
1,\ldots,n\right\}  ,$ $\left[  n..m\right]  =\left\{  n,n+1,\ldots,m\right\}
.$ Write $X^{\left(  k\right)  }$ for the collection of $k$-sets of a set $X.$

Given a graph $G$ and disjoint nonempty sets $X,Y\subset V\left(  G\right)  ,$
we denote the number of $X-Y$ edges by $e_{G}\left(  X,Y\right)  $ and set
$\sigma_{G}\left(  X,Y\right)  =e_{G}\left(  X,Y\right)  /\left(  \left\vert
X\right\vert \left\vert Y\right\vert \right)  .$ Likewise, $e_{G}\left(
X\right)  $ is the number of edges induced by $X$ and $\sigma_{G}\left(
X\right)  =2e_{G}\left(  X\right)  /\left\vert X\right\vert ^{2}.$

Furthermore, $G\left[  X\right]  $ stands for the graph induced by $X,$
$\Gamma_{G}\left(  X\right)  $ is the set of vertices joined to all $u\in X,$
and $d_{G}\left(  X\right)  =\left\vert \Gamma_{G}\left(  X\right)
\right\vert .$ In any of the functions $e_{G}\left(  X,Y\right)  ,$
$\sigma_{G}\left(  X,Y\right)  ,$ $\sigma_{G}\left(  X\right)  ,$
$e_{G}\left(  X\right)  ,$ $\Gamma_{G}\left(  X\right)  ,$ and $d_{G}\left(
X\right)  $ we drop the subscript if the graph $G$ is understood.

As usual, $\delta\left(  G\right)  $ and $\Delta\left(  G\right)  $ denote the
minimum and maximum degrees of $G$, and $\omega\left(  G\right)  $ denotes its
clique number. We write $\psi\left(  G\right)  $ for the order of the largest
component of $G.$

Given $\varepsilon>0,$ a pair $\left(  A,B\right)  $ of nonempty disjoint sets
$A,B\subset V\left(  G\right)  $ is called $\varepsilon$\emph{-regular} if
$\left\vert \sigma\left(  A,B\right)  -\sigma\left(  X,Y\right)  \right\vert
<\varepsilon$ whenever $X\subset A,$ $Y\subset B,$ $\left\vert X\right\vert
\geq\varepsilon\left\vert A\right\vert ,$ $\left\vert Y\right\vert
\geq\varepsilon\left\vert B\right\vert .$ Given $\varepsilon>0,$ a partition
$V\left(  G\right)  =\cup_{i=0}^{k}V_{i}$ is called $\varepsilon
$\emph{-regular,} if $\left\vert V_{0}\right\vert <\varepsilon n,$ $\left\vert
V_{1}\right\vert =\cdots=\left\vert V_{k}\right\vert ,$ and for every
$i\in\left[  k\right]  ,$ all least $\left(  1-\varepsilon\right)  k$ pairs
$\left(  V_{i},V_{j}\right)  $ are $\varepsilon$-regular for $j\in\left[
k\right]  \backslash\left\{  i\right\}  .$

Let $y,x_{1},\ldots,x_{k}$ be real variables. The notation $y\ll\left(
x_{1},\ldots,x_{k}\right)  $ is equivalent to \textquotedblleft$y>0$ and $y$
is sufficiently small, given $x_{1},\ldots,x_{k}$\textquotedblright\ or, in
other words, \textquotedblleft there exists a function $y_{0}\left(
x_{1},\ldots,x_{k}\right)  >0$ and $0<y\leq y_{0}\left(  x_{1},\ldots
,x_{k}\right)  $\textquotedblright. Likewise, $y\gg\left(  x_{1},\ldots
,x_{k}\right)  $ is equivalent to \textquotedblleft$y$ is sufficiently large,
given $x_{1},\ldots,x_{k}$\textquotedblright\ or, in other words,
\textquotedblleft there exists a function $y_{0}\left(  x_{1},\ldots
,x_{k}\right)  >0$ and $y\geq y_{0}\left(  x_{1},\ldots,x_{k}\right)
$\textquotedblright$.$ Since explicit bounds on $y_{0}\left(  x_{1}%
,\ldots,x_{k}\right)  $ are often cumbersome and of little use, we believe
that the above notation simplifies the presentation and emphasizes the
dependence between the relevant variables.

The following structural result proved in \cite{Nik06} is a key ingredient of
our proof of Theorem \ref{mth}.

\begin{fact}
\label{densth} For all $0<\varepsilon<1,$ $p\geq2,$ there exist $\varsigma
=\varsigma\left(  \varepsilon,p\right)  >0$ and $L=L\left(  \varepsilon
,p\right)  $ such that for every graph $G$ of sufficiently large order $n$
with $k_{p}\left(  G\right)  <\varsigma n^{p},$ there exists a partition
$V\left(  G\right)  =\cup_{i=0}^{L}V_{i}$ with the following properties:

- $\left\vert V_{0}\right\vert <\varepsilon n,$ $\left\vert V_{1}\right\vert
=\cdots=\left\vert V_{L}\right\vert ;$

- $\Delta\left(  G\left[  V_{i}\right]  \right)  <\varepsilon\left\vert
V_{i}\right\vert $ for every $i\in\left[  k\right]  .$
\end{fact}

For a general introduction to the Regularity Lemma of Szemer\'{e}di
\cite{Sze78} the reader is referred to \cite{Bol98} and \cite{KoSi96}. We
shall use the following specific form implied by Fact \ref{densth}.

\begin{fact}
\label{maint} For all $0<\varepsilon<1,$ $p\geq2,$ and $k_{0}\geq2,$ there
exist $\rho=\rho\left(  \varepsilon,p,k_{0}\right)  >0$ and $K=K\left(
\varepsilon,p,k_{0}\right)  $ such that for every graph $G$ of sufficiently
large order $n$ with $k_{p+1}\left(  G\right)  <\rho n^{p+1}$, there exists an
$\varepsilon$-regular partition $V\left(  G\right)  =\cup_{i=0}^{k}V_{i}$ with
$k_{0}\leq k<K,$ and $\Delta\left(  G\left[  V_{i}\right]  \right)
<\varepsilon\left\vert V_{i}\right\vert $ for every $i\in\left[  k\right]  .$
\end{fact}

Also we shall use the following simplified versions of the Counting Lemma.

\begin{fact}
\label{tKL} Let $0<\varepsilon<d<1$ and $\left(  A,B\right)  $ be an
$\varepsilon$-regular pair with $\sigma\left(  A,B\right)  \geq d.$ Then there
are at least $\left(  1-\varepsilon\right)  \left\vert A\right\vert $ vertices
$v\in A$ with $\left\vert \Gamma\left(  v\right)  \cap B\right\vert \geq
d-\varepsilon$.
\end{fact}

\bigskip

\begin{fact}
\label{ctL} For all $0<d<1$ and $p\geq2,$ there exist $\varepsilon_{0}$ and
$t_{0}$ such that the following assertion holds:

Let $\varepsilon>\varepsilon_{0},$ $t>t_{0},$ $G$ be a graph of order $pt,$
and $V(G)=\cup_{i=1}^{p}V_{i}$ be a partition such that $\left\vert
V_{1}\right\vert =\cdots=\left\vert V_{p}\right\vert =t.$ If for every $1\leq
i<j\leq p$ the pair $(V_{i},V_{j})$ is $\varepsilon$-regular and $\sigma
(V_{i},V_{j})\geq d,$ then $k_{p}\left(  G\right)  \geq d^{p^{2}}t^{p}.$
\end{fact}

The following lemma can be traced back to Kostochka and R\"{o}dl
\cite{KoRo01}; it was used later by other researchers in various forms, see,
e.g., \cite{Sud05} and its references. We prove the lemma in \ref{pprobl}.

\begin{lemma}
\label{probl} For all $k\geq2,$ $d>0,$ $\lambda>0$ there exists $a=a\left(
k,d,\lambda\right)  >0$ such that for every graph $G$ and nonempty disjoint
sets $U_{1},U_{2}\subset V\left(  G\right)  $ with $e\left(  G\right)  \geq
d\left\vert U_{1}\right\vert \left\vert U_{2}\right\vert ,$ and sufficiently
large $\left\vert U_{1}\right\vert $ there exists $W\subset U_{1}$ with
$\left\vert W\right\vert \geq\left\vert U_{1}\right\vert ^{1-\lambda}$ and
$d\left(  X\right)  >a\left\vert U_{2}\right\vert $ for every $X\subset
W^{\left(  k\right)  }.$
\end{lemma}

\subsection{Proof of Theorem \ref{mth}}

Set $N=p\left(  n-1\right)  +1$ and let $E\left(  K_{N}\right)  =E\left(
R\right)  \cup E\left(  B\right)  $ be a $2$-coloring. In the following list
we show how the variables used in our proof depend on each other
\begin{align*}
\alpha &  \ll\left(  p,q\right)  ,\\
\theta &  \ll\left(  p,q,\alpha\right)  ,\\
\xi &  \ll\left(  p,q,\alpha,\gamma\right)  ,\\
\beta &  \ll\left(  p,q,\alpha,\gamma,\xi\right)  ,\\
\varepsilon &  \ll\left(  p,q,\alpha,\beta,\gamma,\xi\right)  ,\\
k_{0}  &  \gg\left(  p,q,\alpha,\beta,\gamma,\xi,\varepsilon\right)  ,\\
\eta &  \ll\left(  p,q,k_{0},\alpha,\beta,\gamma,\xi,\varepsilon\right)  ,\\
c  &  \ll\left(  p,q,k_{0},\alpha,\beta,\gamma,\xi,\varepsilon\right)  ,\\
n  &  \gg\left(  p,q,k_{0},\alpha,\beta,\gamma,\xi,\varepsilon\right)  .
\end{align*}

Let $K\left(  \varepsilon,p,k_{0}\right)  ,$ $\rho\left(  \varepsilon
,p,k_{0}\right)  ,$ $\varsigma\left(  \varepsilon,p\right)  ,$ $L\left(
\varepsilon,p\right)  ,$ $a\left(  k,d,\lambda\right)  $ be as defined in Fact
\ref{densth}, Fact \ref{maint}, and Lemma \ref{probl}. Assume that
\begin{equation}
js_{p+1}\left(  R\right)  \leq cn^{p-1} \label{massum}%
\end{equation}
and select a $q$-degenerate $\left(  \gamma,\eta\right)  $-splittable graph
$H$ of order $n.$ To prove the theorem, we shall show that $H\subset B$.

Assumption (\ref{massum}) implies that
\[
k_{p+1}\left(  R\right)  \leq\binom{p+1}{2}^{-1}\binom{N}{2}js_{p+1}\left(
R\right)  <cN^{p+1}<\rho\left(  \varepsilon,p+1,k_{0}\right)  N^{p+1}.
\]
Thus, by Fact \ref{maint}, there exists an $\varepsilon$-regular partition
$V\left(  R\right)  =\cup_{i=0}^{k}V_{i}$ such that $k_{0}<k<K\left(
\varepsilon,p+1,k_{0}\right)  $ and $\Delta\left(  R\left[  V_{i}\right]
\right)  <\varepsilon\left\vert V_{i}\right\vert $ for all $i\in\left[
k\right]  .$ Set $t=\left\vert V_{1}\right\vert =\cdots=\left\vert
V_{k}\right\vert $ and note that for all $i\in\left[  k\right]  ,$
\begin{equation}
\delta\left(  B\left[  V_{i}\right]  \right)  =t-1-\Delta\left(  R\left[
V_{i}\right]  \right)  >t-\varepsilon t-1>\left(  1-\beta\right)  t.
\label{minB}%
\end{equation}

Next define the graphs $R^{\ast}$ and $B^{\ast}$ by%
\begin{align*}
V\left(  B^{\ast}\right)   &  =\left[  k\right]  ,\text{ \ \ }V\left(
R^{\ast}\right)  =\left[  k\right]  ,\text{ }\\
E\left(  B^{\ast}\right)   &  =\left\{  \left\{  u,v\right\}  :1\leq u<v\leq
k\text{ and }\sigma_{B}\left(  V_{u},V_{v}\right)  >1-\beta\right\}  ,\\
E\left(  R^{\ast}\right)   &  =\left\{  \left\{  u,v\right\}  :1\leq u<v\leq
k,\text{ }\left(  V_{u},V_{v}\right)  \text{ is }\varepsilon\text{-regular and
}\sigma_{R}\left(  V_{u},V_{v}\right)  \geq\beta\right\}  .
\end{align*}

Note first that $E\left(  B^{\ast}\right)  \cap E\left(  R^{\ast}\right)
=\varnothing$. Moreover, for every vertex $u\in\left[  k\right]  ,$ we have
\begin{equation}
d_{B^{\ast}}\left(  u\right)  +d_{R^{\ast}}\left(  u\right)  >k-1-\varepsilon
k. \label{md}%
\end{equation}
Indeed, if $\left\{  u,v\right\}  \notin E\left(  B^{\ast}\right)  \cup
E\left(  R^{\ast}\right)  $ then the pair $\left(  V_{u},V_{v}\right)  $ is
not $\varepsilon$-regular; hence $\left\{  u,v\right\}  \notin E\left(
B^{\ast}\right)  \cup E\left(  R^{\ast}\right)  $ holds for fewer than
$\varepsilon k$ vertices $v\in\left[  k\right]  \backslash\left\{  u\right\}
$.

We first show that $H\subset B$ if $\Delta\left(  B^{\ast}\right)  $
satisfies
\begin{equation}
\Delta\left(  B^{\ast}\right)  \geq\left(  1+2\xi\right)  \frac{k}{p}.
\label{maxd0}%
\end{equation}
Indeed, set $r=\Delta\left(  B^{\ast}\right)  $ and select $v_{0}\in\left[
k\right]  $ with $d_{B^{\ast}}\left(  v_{0}\right)  =r.$ Let $\Gamma_{B^{\ast
}}\left(  v_{0}\right)  =\left\{  v_{1},\ldots,v_{r}\right\}  $ and set
$U_{j}=V_{v_{j}}$ for $j=0,\ldots,r.$

To simplify the presentation of our proof we formulate various claims proved
later in \ref{Claims}.

\begin{claim}
\label{cl1}$H\subset B\left[  \cup_{i=0}^{r}U_{i}\right]  .$
\end{claim}

Hereafter we shall assume that (\ref{maxd0}) fails, i.e.,%
\begin{equation}
\Delta\left(  B^{\ast}\right)  <\left(  1+2\xi\right)  \frac{k}{p}.
\label{maxd1}%
\end{equation}
In view of (\ref{md}), this inequality implies a lower bound on $\delta\left(
R^{\ast}\right)  ,$ viz.%
\begin{equation}
\delta\left(  R^{\ast}\right)  >k-1-\varepsilon k-\left(  1+2\xi\right)
\frac{k}{p}>\left(  \frac{p-1}{p}-2\xi\right)  k. \label{mindR}%
\end{equation}
In turn, the bound (\ref{mindR}), together with the assumption (\ref{massum}),
implies a definite structure in $R^{\ast}.$

\begin{claim}
\label{cl2}$R^{\ast}$ is $p$-partite.
\end{claim}

Write $Z_{1},\ldots,Z_{p}$ for the color classes of $R^{\ast}.$ For every
$i\in\left[  k\right]  ,$ let $\mu\left(  i\right)  \in\left[  p\right]  $ be
the unique value satisfying $i\in Z_{\mu\left(  i\right)  }.$ Observe that the
sets $Z_{1},\ldots,Z_{p}$ determine a partition of $\left[  N\right]
\backslash V_{0}$ into $p$ sets that are dense in $B.$ Indeed, $e_{B}\left(
V_{u}\right)  >\left(  1-\beta\right)  t^{2}/2$ for all $u\in\left[  k\right]
$, and also $e_{B}\left(  V_{u},V_{v}\right)  >\left(  1-\beta\right)  t^{2}$
whenever $\mu\left(  u\right)  =\mu\left(  v\right)  $ and $u\neq v.$

Next we show that the color classes of $R^{\ast}$ cannot be two small. Indeed,
in view of (\ref{maxd1}), for every $i\in\left[  p\right]  ,$ we have%
\begin{align}
\left\vert Z_{i}\right\vert  &  =k-%
{\textstyle\sum\limits_{j\in\left[  p-1\right]  \backslash\left\{  i\right\}
}}
\left\vert Z_{j}\right\vert \geq k-\left(  p-1\right)  \left(  \Delta\left(
B^{\ast}\right)  +1\right) \label{minz}\\
&  =k-\left(  1+2\xi\right)  \frac{\left(  p-1\right)  k}{p}-p+1>\left(
1-2p\xi\right)  \frac{k}{p}.\nonumber
\end{align}

In Claims \ref{cl3}-\ref{cl6} we show that $H\subset B$ provided the
inequality
\begin{equation}%
{\displaystyle\sum\limits_{1\leq h<s\leq p}}
\left(
{\displaystyle\sum\limits_{i\in Z_{h},\text{ }j\in Z_{s}}}
e_{B}\left(  V_{i},V_{j}\right)  \right)  \geq\frac{\alpha}{2}N^{2}
\label{assum2}%
\end{equation}
holds. Inequality (\ref{assum2}) implies that $e_{B}\left(  V_{u}%
,V_{v}\right)  $ is substantial for substantially many pairs $u,v$ with
$\mu\left(  u\right)  \neq\mu\left(  v\right)  ;$ we shall use this fact to
embed a substantial part of $H.$ Let us first derive a more specific condition
from (\ref{assum2}).

\begin{claim}
\label{cl3}There exist $i_{1}\in\left[  k\right]  $ and $j\in\left[  p\right]
\backslash\mu\left(  i_{1}\right)  $ such that
\[
\sum_{\mu\left(  s\right)  =j}e_{B}\left(  V_{i_{1}},V_{s}\right)
>\alpha\left\vert Z_{j}\right\vert t^{2}.
\]

\end{claim}

We may and shall assume that $i_{1}\in Z_{1}$ and $j=2.$ Setting
$X=\cup\left\{  V_{s}:s\in Z_{2}\right\}  ,$ we see that Claim \ref{cl3}
amounts to
\begin{equation}
e_{B}\left(  V_{i_{1}},X\right)  >\alpha\left\vert Z_{2}\right\vert t^{2}.
\label{cl1in}%
\end{equation}

Observe also that, in view of (\ref{minz}), we have
\begin{equation}
\left\vert X\right\vert =\left\vert Z_{2}\right\vert t\geq\left(
1-2p\xi\right)  \frac{kt}{p}. \label{minx}%
\end{equation}
In addition,
\begin{align*}
2e_{B}\left(  X\right)   &  =2\sum_{s\in Z_{2}}e_{B}\left(  V_{s}\right)
+\sum_{i\in Z_{2}}\left(  \sum_{j\in Z_{2}\backslash\left\{  i\right\}  }%
e_{B}\left(  V_{i},V_{j}\right)  \right) \\
&  >\left\vert Z_{2}\right\vert \left(  1-\beta\right)  t^{2}+\left\vert
Z_{2}\right\vert \left(  \left\vert Z_{2}\right\vert -1\right)  \left(
1-\beta\right)  t^{2}=\left\vert Z_{2}\right\vert ^{2}\left(  1-\beta\right)
t^{2},
\end{align*}
and so%
\begin{equation}
\sigma_{B}\left(  X\right)  >\left(  1-\beta\right)  . \label{minex}%
\end{equation}

Inequality (\ref{cl1in}) implies that substantially many vertices in
$V_{i_{1}}$ are joined to substantially many vertices in $X.$ In the following
claim we strengthen this condition.

\begin{claim}
\label{cl3.0} There exists $W_{0}\subset V_{i_{1}}$ with $\left\vert
W_{0}\right\vert >\left(  \alpha/2\right)  t$ such that for all $u\in W_{0},$
\[
\left\vert \Gamma_{B}\left(  u\right)  \cap X\right\vert >\left(
\alpha/2\right)  \left\vert X\right\vert .
\]

\end{claim}

Next set $Y=\cup\left\{  V_{s}:s\in Z_{1},\text{ }s\neq i_{1}\right\}  ;$ by
(\ref{minz}) we have%
\begin{equation}
\left\vert Y\right\vert =\left(  \left\vert Z_{1}\right\vert -1\right)
t\geq\left\vert Z_{1}\right\vert \left(  1-\frac{1}{\left\vert Z_{1}%
\right\vert }\right)  t\geq\left\vert Z_{1}\right\vert \left(  1-\frac
{p}{\left(  1+2\xi\right)  k_{0}}\right)  t\geq\left(  1-\beta\right)
\left\vert Z_{1}\right\vert t. \label{miny}%
\end{equation}

In addition,
\begin{align*}
2e_{B}\left(  Y\right)   &  =2\sum_{s\in Z_{1}\backslash\left\{
i_{1}\right\}  }e_{B}\left(  V_{s}\right)  +\sum_{i\in Z_{1}\backslash\left\{
i_{1}\right\}  }\left(  \sum_{j\in Z_{1}\backslash\left\{  i,i_{1}\right\}
}e_{B}\left(  V_{i},V_{j}\right)  \right) \\
&  >\left(  \left\vert Z_{1}\right\vert -1\right)  \left(  1-\beta\right)
t^{2}+\left(  \left\vert Z_{1}\right\vert -1\right)  \left(  \left\vert
Z_{1}\right\vert -2\right)  \left(  1-\beta\right)  t^{2}=\left(  \left\vert
Z_{1}\right\vert -1\right)  ^{2}\left(  1-\beta\right)  t^{2},
\end{align*}
and so%
\begin{equation}
\sigma_{B}\left(  Y\right)  >\left(  1-\beta\right)  . \label{miney}%
\end{equation}

Inequality (\ref{miny}) implies that substantially many vertices in $W_{0}$
are joined to substantially many vertices in $Y.$ Next we strengthen this condition.

\begin{claim}
\label{cl3.1} There exists $W_{1}\subset W_{0}$ with $\left\vert
W_{1}\right\vert >\left(  \alpha/4\right)  t$ such that for all $u\in W_{1},$
\[
\left\vert \Gamma_{B}\left(  u\right)  \cap Y\right\vert >\left(
1-\sqrt{\beta}\right)  \left\vert Y\right\vert .
\]

\end{claim}

Furthermore, the lower bound on $\delta\left(  R^{\ast}\right)  $ given by
inequality (\ref{mindR}) implies that $i_{1}$ belongs to a $p$-clique in
$R^{\ast}$.

\begin{claim}
\label{cl4}There exist $i_{2}\in Z_{2},\ldots,i_{p}\in Z_{p}$ such that
$\left\{  i_{1},i_{2},\ldots,i_{p}\right\}  $ induces a clique in $R^{\ast}.$
\end{claim}

Claim \ref{cl4}, together with $js_{p+1}\left(  R\right)  <cn^{p-1},$ implies
that the graph $B\left[  W_{1}\right]  $ contains a large clique.

\begin{claim}
\label{cl5}There exists $W\subset W_{1}$ with $\left\vert W\right\vert \geq
t^{1-\gamma/2}$ such that $B\left[  W\right]  $ is a complete graph.
\end{claim}

In summary, Claims \ref{cl3}-\ref{cl5} together with (\ref{minx}) and
(\ref{miny}) imply that the sets $W,$ $X,$ and $Y$ have the following properties:

- $\left\vert W\right\vert \geq t^{1-\gamma/2}$ and $B\left[  W\right]  $ is a
complete graph,

- $\left\vert X\right\vert \geq\left(  1-2p\xi\right)  k/p,$

- $\left\vert Y\right\vert \geq\left(  1-2p\xi\right)  k/p,$

- $\left\vert \Gamma_{B}\left(  u\right)  \cap X\right\vert >\left(
\alpha/4\right)  \left\vert X\right\vert $ and $\left\vert \Gamma_{B}\left(
u\right)  \cap Y\right\vert >\left(  1-\sqrt{\beta}\right)  \left\vert
Y\right\vert ,$ for every $u\in W.$

It turns out that these properties are sufficient to achieve our goal - to
embed $H$.

\begin{claim}
\label{cl6}$H\subset B\left[  W\cup X\cup Y\right]  .$
\end{claim}

Hereafter we shall assume that (\ref{assum2}) fails, i.e.,
\begin{equation}
\sum_{1\leq h<s\leq p}\left(  \sum_{i\in Z_{h},\text{ }j\in Z_{s}}e_{B}\left(
V_{i},V_{j}\right)  \right)  <\frac{\alpha}{2}N^{2}. \label{loblu}%
\end{equation}

This inequality implies that $e_{B}\left(  V_{u},V_{v}\right)  $ is small for
most pairs $u,v$ with $\mu\left(  u\right)  \neq\mu\left(  v\right)  .$ We
shall deduce that $R$ can be made $p$-partite by removing only a small
proportion of its vertices.

\begin{claim}
\label{cl7} $R$ contains an induced $p$-partite subgraph $R_{1}$ with color
classes $U_{1},\ldots,U_{p}$ such that $\left\vert U_{1}\right\vert
=\cdots=\left\vert U_{p}\right\vert >\left(  1-\theta\right)  n$ and
\[
\left\vert \Gamma_{R_{1}}\left(  u\right)  \cap U_{i}\right\vert >\left(
1-\theta\right)  n
\]
\ for each $i\in\left[  p\right]  $ and $u\in V\left(  R_{1}\right)
\backslash U_{i}.$
\end{claim}

Since $R_{1}$ is an induced $p$-partite subgraph of $R,$ the graph $B$
contains cliques of size close to $n;$ hence $H$ can be embedded in $B$ almost
entirely; to embed $H$ in full, we need an additional argument. Analyzing the
way vertices from $V\left(  R\right)  \backslash V\left(  R_{1}\right)  $ can
be joined to the vertices of $R_{1},$ we derive the following assertion.

\begin{claim}
\label{cl8}There exist disjoint sets $M\subset V\left(  R_{1}\right)  $ and
$A,C\subset V\left(  R\right)  \backslash V\left(  R_{1}\right)  $ such that%
\begin{align}
\left\vert M\right\vert +\left\vert A\right\vert +\left\vert C\right\vert  &
=n-1+\left\lceil \frac{\left\vert C\right\vert +1}{p}\right\rceil
,\label{cond}\\
\left\vert A\right\vert  &  <\theta n,\label{cond2}\\
\left\vert C\right\vert  &  <2\theta n \label{cond3}%
\end{align}
with the following properties:

(i) $B\left[  M\right]  $ is a complete graph;

(ii) $\Gamma_{B}\left(  u\right)  \cap M=M$ for every vertex $u\in A;$

(iii) $\left\vert \Gamma_{B}\left(  u\right)  \cap M\right\vert \geq\left(
1-p^{2}\theta\right)  \left\vert M\right\vert $ for every vertex $u\in C.$
\end{claim}

Using the properties of the sets $M,A,$ and $C$ we embed $H$, completing the
proof of the theorem.

\begin{claim}
\label{cl9} $H\subset B\left[  M\cup A\cup C\right]  $.
\end{claim}

\subsubsection{Results supporting proofs of the claims}

\bigskip

\begin{fact}
Every subgraph of a $q$-degenerate graph is $q$-degenerate.
\end{fact}

\begin{fact}
\label{arang}The vertices of any $q$-degenerate graph of order $n$ can be
labeled $\left\{  v_{1},\ldots,v_{n}\right\}  $ so that $\left\vert
\Gamma\left(  v_{i}\right)  \cap\left\{  v_{1},\ldots,v_{i-1}\right\}
\right\vert \leq q$ for every $i\in\left[  n\right]  .$
\end{fact}

\begin{fact}
\label{color}Every $q$-degenerate graph is $\left(  q+1\right)  $-partite.
\end{fact}

\begin{proposition}
\label{cheb}In any $q$-degenerate graph $H$ the number of vertices of degree
$2q+1$ or higher is at most $2q\left\vert H\right\vert /\left(  2q+1\right)
.$
\end{proposition}

\begin{proof}
Letting $S=\left\{  u:u\in V\left(  H\right)  ,\text{ }d\left(  u\right)
\geq2q+1\right\}  ,$ we have
\[
2q\left\vert H\right\vert \geq2e\left(  H\right)  \geq\sum_{u\in V\left(
H\right)  }d\left(  u\right)  \geq\sum_{u\in S}d\left(  u\right)  \geq\left(
2q+1\right)  \left\vert S\right\vert ,
\]
and the assertion follows.
\end{proof}

\bigskip

\begin{lemma}
\label{embedq}Let $q\geq0,$ $\tau>0,$ and $G=G\left(  n\right)  $ be a graph
with $\delta\left(  G\right)  \geq\left(  1-\tau\right)  n.$ Then $G$ contains
all $q$-degenerate graphs of order $l\leq\left(  1-q\tau\right)  n.$
\end{lemma}

\begin{proof}
We use induction on $l.$ The assertion holds trivially for $l=1;$ assume that
it holds for $1\leq l^{\prime}<l.$ Let $H$ be a $q$-degenerate graph of order
$l$ and $u\in V\left(  H\right)  $ be a vertex with $d_{H}\left(  u\right)
=d\leq q.$ Let $\Gamma_{H}\left(  u\right)  =\left\{  v_{1},\ldots
,v_{d}\right\}  $ and $H^{\prime}=H-u.$ By the induction assumption there
exists a monomorphism $\varphi:H^{\prime}\rightarrow G.$ We have
\begin{align*}
\left\vert
{\textstyle\bigcap_{i=1}^{d}}
\Gamma_{G}\left(  \varphi\left(  v_{i}\right)  \right)  \right\vert  &
\geq\sum_{i=1}^{d}d_{G}\left(  \varphi\left(  v_{i}\right)  \right)  -\left(
d-1\right)  n>d\left(  1-\tau\right)  n-\left(  d-1\right)  n\\
&  \geq\left(  1-q\tau\right)  n>l^{\prime}.
\end{align*}
Hence there exists $v\in\left(  \bigcap_{i=1}^{d}\Gamma_{G}\left(
\varphi\left(  v_{i}\right)  \right)  \right)  \backslash\varphi\left(
V\left(  H^{\prime}\right)  \right)  .$ To complete the induction step and the
proof, define a monomorphism $\varphi^{\prime}:H\rightarrow G$ by
\[
\varphi^{\prime}\left(  w\right)  =\left\{
\begin{array}
[c]{ll}%
\varphi\left(  w\right)  , & \text{if }w\in V\left(  H^{\prime}\right) \\
v, & \text{if }w=u.
\end{array}
\right.
\]

\end{proof}

\bigskip

\begin{lemma}
\label{triml}Suppose $G$ is a $\left(  \gamma,\eta\right)  $-splittable
$q$-degenerate graph of order $n$. Then there exists $M\subset V\left(
G\right)  $ such that $\left\vert M\right\vert <\left(  2q+1\right)
n^{1-\gamma},$ and $\psi\left(  G-M\right)  <\eta n$ and $\left\vert
\Gamma\left(  u\right)  \cap M\right\vert \leq2q$ for every $u\in V\left(
G\right)  \backslash M.$
\end{lemma}

\begin{proof}
Since $G$ is $\left(  \gamma,\eta\right)  $-splittable, there is a set
$N\subset V\left(  G\right)  $ such that $\left\vert N\right\vert
<n^{1-\gamma}$ and $\psi\left(  G-N\right)  <\eta n.$ Set $M=N$ and apply the
following procedure to $G:$

\emph{While there exists }$u\in V\left(  G\right)  \backslash M$\emph{ with
}$\left\vert \Gamma\left(  u\right)  \cap M\right\vert \geq2q+1$ \emph{do }

$\qquad M:=M\cup\left\{  u\right\}  $

\emph{end.}

Set $M^{\prime}=\left\{  u:u\in M,\text{ }\left\vert \Gamma\left(  u\right)
\cap M\right\vert \geq2q+1\right\}  .$ Since $G\left[  M\right]  $ is
$q$-degenerate, Proposition \ref{cheb} implies that $\left\vert M^{\prime
}\right\vert \leq2q\left\vert M\right\vert /\left(  2q+1\right)  .$ By our
selection, $\left\vert \Gamma\left(  u\right)  \cap M\right\vert \geq2q+1$ for
all of $u\in M\backslash N;$ hence, $\left\vert M\backslash N\right\vert
\leq2q\left\vert M\right\vert /\left(  2q+1\right)  ,$ implying that
$\left\vert M\right\vert \leq\left(  2q+1\right)  \left\vert N\right\vert
\leq\left(  2q+1\right)  n^{1-\gamma}.$
\end{proof}

\bigskip

\begin{proposition}
\label{dendeg} Let $0<\tau<1$ and $G$ be a graph of order $n$ with $e\left(
G\right)  >\left(  1-\tau\right)  n^{2}/2.$ Then $G$ contains an induced
subgraph $G_{0}$ with $\left\vert G_{0}\right\vert >\left(  1-\sqrt{\tau
}\right)  n$ and $\delta\left(  G_{0}\right)  >\left(  1-2\sqrt{\tau}\right)
n$
\end{proposition}

\begin{proof}
Let
\[
W=\left\{  u:d_{G}\left(  u\right)  >\left(  1-\sqrt{\tau}\right)  n\right\}
.
\]
We have
\begin{align*}
\left(  1-\tau\right)  n^{2}  &  <2e\left(  G\right)  =\sum_{u\in W}%
d_{G}\left(  u\right)  +\sum_{u\in V\left(  G\right)  \backslash W}%
d_{G}\left(  u\right)  \leq n\left\vert W\right\vert +\left(  1-\sqrt{\tau
}\right)  n\left(  n-\left\vert W\right\vert \right) \\
&  =\sqrt{\tau}n\left\vert W\right\vert +\left(  1-\sqrt{\tau}\right)  n^{2},
\end{align*}
and so $\left(  1-\sqrt{\tau}\right)  n<\left\vert W\right\vert .$
Furthermore, for every $u\in W,$
\[
\left\vert \Gamma_{G}\left(  u\right)  \cap W\right\vert \geq\left\vert
\Gamma_{G}\left(  u\right)  \right\vert -\left\vert V\left(  G\right)
\backslash W\right\vert \geq\left(  1-2\sqrt{\tau}\right)  n.
\]
Thus, setting $G_{0}=G\left[  W\right]  ,$ the proof is completed.
\end{proof}

\bigskip

\begin{fact}
[\cite{BoNi04}]\label{stabj} Let $p\geq3,$ $n>p^{8}$, and $0<\alpha
<p^{-8}/16.$ If a graph $G=G\left(  n\right)  $ satisfies
\[
e\left(  G\right)  >\left(  \frac{p-1}{2p}-\alpha\right)  n^{2},
\]
then either
\begin{equation}
js_{p}\left(  G\right)  >\left(  1-\frac{1}{p^{3}}\right)  \frac{n^{p-2}%
}{p^{p+5}}, \label{minjs}%
\end{equation}
or $G$ contains an induced $p$-partite subgraph $G_{0}$ of order at least
$\left(  1-2\sqrt{\alpha}\right)  n$ with minimum degree
\begin{equation}
\delta\left(  G_{0}\right)  >\left(  1-\frac{1}{p}-4\sqrt{\alpha}\right)  n.
\label{mindg}%
\end{equation}

\end{fact}

\bigskip

\begin{fact}
[\cite{BoNi04}]\label{minj}Let $2\leq r<\omega\left(  G\right)  $ and
$\alpha\geq0.$ If $G=G\left(  n\right)  $ and%
\[
\delta\left(  G\right)  \geq\left(  \frac{r-1}{r}+\alpha\right)  n
\]
then%
\[
k_{r+1}\left(  G\right)  \geq\alpha\frac{r^{2}}{r+1}\left(  \frac{n}%
{r}\right)  ^{r+1}.
\]

\end{fact}

\subsubsection{\label{Claims}Proofs of the claims}

Let $K\left(  \varepsilon,p,k\right)  $ and $\rho\left(  \varepsilon
,p,k\right)  ,$ $\varsigma\left(  \varepsilon,p\right)  ,$ and $L\left(
\varepsilon,p\right)  ,$ be as defined in Fact \ref{densth} and Fact
\ref{maint}; set $K=K\left(  \varepsilon,p,k_{0}\right)  .$

\begin{proof}
[\textbf{Proof of Claim \ref{cl1}}]Set $\varsigma=\varsigma\left(  1/\left(
2q\right)  ,p\right)  $ and $L=L\left(  1/\left(  2q\right)  ,p\right)  .$

Note first that the sets $U_{0},\ldots,U_{r}$ satisfy the following conditions:

- $\left\vert U_{0}\right\vert =\ldots=\left\vert U_{r}\right\vert =t$;

-\emph{ }$\delta\left(  B\left[  U_{i}\right]  \right)  \geq\left(
1-\varepsilon\right)  t>\left(  1-\beta\right)  t$ for $i=0,\ldots,r;$

- $\sigma_{B}\left(  U_{0},U_{i}\right)  >1-\beta$ for $i=1,\ldots,r.$

For every $u\in U_{0}$ set
\[
D\left(  u\right)  =\left\{  i:i\in\left[  r\right]  ,\text{ }\left\vert
\Gamma\left(  u\right)  \cap U_{i}\right\vert \geq\left(  1-2\sqrt{\beta
}\right)  t\right\}
\]
and let
\[
W=\left\{  u:u\in U_{0},\text{ }\left\vert D\left(  u\right)  \right\vert
\geq\left(  1-\sqrt{\beta}\right)  r\right\}  .
\]
We shall prove that $\left\vert W\right\vert >t/2.$ Indeed, we see that
\begin{align*}
\left(  1-\beta\right)  rt^{2}  &  <\sum_{i\in\left[  r\right]  }e\left(
U_{0},U_{i}\right)  =\sum_{u\in U_{0}}\left(  \sum_{i\in\left[  r\right]
}\left\vert \Gamma\left(  u\right)  \cap U_{i}\right\vert \right) \\
&  =\sum_{u\in W}\left(  \sum_{i\in\left[  r\right]  }\left\vert \Gamma\left(
u\right)  \cap U_{i}\right\vert \right)  +\sum_{u\in U_{0}\backslash W}\left(
\sum_{i\in\left[  r\right]  }\left\vert \Gamma\left(  u\right)  \cap
U_{i}\right\vert \right) \\
&  <\left\vert W\right\vert rt+\sum_{u\in U_{0}\backslash W}\left(  D\left(
u\right)  t+\left(  r-D\left(  u\right)  \right)  \left(  1-2\sqrt{\beta
}\right)  t\right) \\
&  \leq\left\vert W\right\vert rt+t\left(  t-\left\vert W\right\vert \right)
\left(  r\left(  1-2\sqrt{\beta}\right)  +2\sqrt{\beta}D\left(  u\right)
\right) \\
&  <\left\vert W\right\vert rt+tr\left(  t-\left\vert W\right\vert \right)
\left(  \left(  1-2\sqrt{\beta}\right)  +2\sqrt{\beta}\left(  1-\sqrt{\beta
}\right)  \right) \\
&  =\left\vert W\right\vert rt+rt\left(  t-\left\vert W\right\vert \right)
\left(  1-2\beta\right)  .
\end{align*}
Hence%
\[
\left(  1-\beta\right)  t\leq\left\vert W\right\vert +\left(  t-\left\vert
W\right\vert \right)  \left(  1-2\beta\right)  =t\left(  1-2\beta\right)
+2\beta\left\vert W\right\vert ,
\]
and so $\left\vert W\right\vert >t/2.$

Since $D\left(  u\right)  \subset\left[  r\right]  ,$ the pigeonhole principle
gives $D\subset\left[  r\right]  $ and $X\subset W$ such that
\[
\left\vert X\right\vert \geq\left\vert W\right\vert /2^{r}\geq t/2^{K+1}%
\]
and $D\left(  u\right)  =D$ for every $u\in X.$ Since%
\begin{align*}
js_{p+1}\left(  R\left[  X\right]  \right)   &  <cn^{p-1}\leq c\left(
\frac{N}{p}\right)  ^{p-1}\leq c\left(  \frac{Kt}{p\left(  1-\varepsilon
\right)  }\right)  ^{p-1}\leq c\left(  Kt\right)  ^{p-1}\\
&  <c\left(  K2^{K+1}\left\vert X\right\vert \right)  ^{p-1}\leq
\varsigma\left\vert X\right\vert ^{p-1},
\end{align*}
Theorem \ref{densth} implies that $X$ contains a set $Y$ with $\left\vert
Y\right\vert \geq\left\vert X\right\vert /2L$ and
\begin{equation}
\delta\left(  B\left[  Y\right]  \right)  >\left(  1-1/2q\right)  \left\vert
Y\right\vert . \label{mind2}%
\end{equation}

On the other hand, Lemma \ref{triml} implies that there exists $M\subset
V\left(  H\right)  $ with $\left\vert M\right\vert \leq\left(  2q+1\right)
\left\vert H\right\vert ^{1-\gamma}$ such that $\psi\left(  H-M\right)
\leq\gamma\left\vert H\right\vert $ and $\left\vert \Gamma_{H}\left(
u\right)  \cap M\right\vert \leq2q$ for every $u\in V\left(  H-M\right)  .$
Since the graph $H\left[  M\right]  $ is $q$-degenerate, we have
\[
\left\vert M\right\vert \leq\left(  2q+1\right)  \left\vert H\right\vert
^{1-\gamma}\leq\left(  2q+1\right)  \left(  rt\right)  ^{1-\gamma}\leq\frac
{t}{2^{r+3}L}\leq\frac{\left\vert X\right\vert }{4L}\leq\frac{\left\vert
Y\right\vert }{2}%
\]
for $t$ large. Hence, in view of (\ref{mind2}), Lemma \ref{embedq} implies
that there exists a monomorphism $\varphi:H\left[  M\right]  \rightarrow Y.$
We shall extend $\varphi$ to $H$ by mapping each component of $H-M$ in turn.

Select a component $C$ of $H-M.$ The choice of the set $M$ implies that%
\[
\left\vert C\right\vert \leq\psi\left(  H-M\right)  \leq\eta\left\vert
H\right\vert <\frac{\sqrt{\beta}\left\vert H\right\vert }{K}<\frac{\sqrt
{\beta}rt}{r}=\sqrt{\beta}t.
\]
We shall extend $\varphi$ over $C$ by mapping $C$ in any set $U_{i},$ $i\in D$
in which there are at least $\left(  6q+1\right)  \sqrt{\beta}t$ vertices
outside of the current range of $\varphi$. Set $l=\left\vert C\right\vert ;$
Proposition \ref{arang} implies that the vertices of $C$ can be arranged as
$v_{1},\ldots,v_{l}$ so that $\left\vert \Gamma_{H}\left(  v_{i}\right)
\cap\left\{  v_{1},\ldots,v_{i-1}\right\}  \right\vert \leq q$ for every
$i\in\left[  l\right]  $. We shall extend $\varphi$ over $C$ mapping each
$v_{i}\in V\left(  C\right)  $ in turn. Suppose we have mapped $v_{1}%
,\ldots,v_{i-1}.$ The vertex $v_{i}$ is joined to at most $q$ vertices from
$\left\{  v_{1},\ldots,v_{i-1}\right\}  $ and at most $2q$ vertices from $M,$
i.e.,
\[
v_{i}\in\left(
{\textstyle\bigcap_{j=1}^{h}}
\Gamma_{H}\left(  v_{i_{j}}\right)  \right)  \cap\left(
{\textstyle\bigcap_{j=1}^{s}}
\Gamma_{H}\left(  u_{i_{j}}\right)  \right)  ,
\]
where $v_{i_{1}},\ldots,v_{i_{h}}\in\left\{  v_{1},\ldots,v_{i-1}\right\}  ,$
$h\leq q,$ and $u_{i_{1}},\ldots,u_{i_{s}}\in M,$ $s\leq2q.$ Set for
convenience $x_{j}=\varphi\left(  v_{i_{j}}\right)  $ for all $j\in\left[
h\right]  ,$ and $y_{j}=\varphi\left(  u_{i_{j}}\right)  $ for all
$j\in\left[  s\right]  .$ Note that
\begin{align*}
&  \left(
{\textstyle\bigcap_{j=1}^{h}}
\Gamma_{B}\left(  x_{j}\right)  \cap U_{i}\right)  \cap\left(
{\textstyle\bigcap_{j=1}^{s}}
\Gamma_{B}\left(  y_{j}\right)  \cap U_{i}\right) \\
&  \geq\sum_{j\in\left[  h\right]  }\left\vert \Gamma_{B}\left(  x_{j}\right)
\cap U_{i}\right\vert +\sum_{j\in\left[  s\right]  }\left\vert \Gamma
_{B}\left(  y_{j}\right)  \cap U_{i}\right\vert -\left(  h+s-1\right)  t\\
&  >\left(  h+s\right)  \left(  1-2\sqrt{\beta}\right)  t-\left(
h+s-1\right)  t\\
&  >\left(  1-6q\sqrt{\beta}\right)  t>\left(  1-\left(  6q+1\right)
\sqrt{\beta}\right)  t+\left\vert C\right\vert .
\end{align*}
Hence there is a vertex $z\in U_{i}$ that is joined to the vertices
$x_{1},\ldots,x_{h},y_{1},\ldots,y_{s}$ and is outside the current range of
$\varphi.$ Setting $\varphi\left(  v_{i}\right)  =z,$ we extend $\varphi$ to a
monomorphism that maps $v_{i}$ into $B$ as well. In this way $\varphi$ can be
extended over the whole component $C$.

Assume for a contradiction that $\varphi$ cannot be extended over some
component $C.$ Therefore, for every $i\in D,$ the current range of $\varphi$
contains at least $\left(  1-\left(  6q+1\right)  \sqrt{\beta}\right)  t$
vertices from $U_{i}.$ Hence%
\begin{align*}
\left\vert H\right\vert  &  \geq\left\vert D\right\vert \left(  1-\left(
6q+1\right)  \sqrt{\beta}\right)  t>\left(  1-\sqrt{\beta}\right)  \left(
1-\left(  6q+1\right)  \sqrt{\beta}\right)  rt\\
&  \geq\left(  1-\left(  6q+2\right)  \sqrt{\beta}\right)  rt>\left(
1-\xi\right)  rt\geq\frac{\left(  1-\xi\right)  \left(  1+2\xi\right)  }%
{p}kt\\
&  \geq\left(  1-\xi\right)  \left(  1+2\xi\right)  \left(  1-\varepsilon
\right)  n>n,
\end{align*}
a contradiction, completing the proof.
\end{proof}

\bigskip

\begin{proof}
[\textbf{Proof of claim \ref{cl2}}]Let $\upsilon=\beta^{p^{2}}.$

We shall prove first that $\omega\left(  R^{\ast}\right)  \leq p.$ Otherwise
by Lemma \ref{ctL} we have
\[
k_{p+1}\left(  R\right)  \geq\upsilon t^{p+1}\geq\upsilon\left(
1-\varepsilon\right)  ^{p+1}\left(  \frac{N}{K}\right)  ^{p+1},
\]
and so
\begin{align*}
js_{p+1}\left(  R\right)   &  \geq\frac{\binom{p+1}{2}}{\binom{N}{2}}%
k_{p+1}\left(  R\right)  >\upsilon\frac{1}{N^{2}}\left(  1-\varepsilon\right)
^{p+1}\left(  \frac{N}{K}\right)  ^{p+1}\\
&  \geq\upsilon\frac{\left(  1-\varepsilon\right)  ^{p+1}}{K^{p+1}}%
N^{p-1}>cn^{p-1},
\end{align*}
contradicting (\ref{massum}).

Since $\omega\left(  R^{\ast}\right)  \leq p,$ and
\[
\delta\left(  R^{\ast}\right)  >\left(  1-\frac{1}{p}-2\xi\right)
k\geq\left(  1-\frac{1}{p-1/3}\right)  k,
\]
by a well-known theorem of Andr\'{a}sfai, Erd\H{o}s, and S\'{o}s \cite{AES74},
$R^{\ast}$ is $p$-partite.
\end{proof}

\bigskip

\begin{proof}
[\textbf{Proof of Claim \ref{cl3}}]In view of (\ref{assum2}), we have%
\begin{align*}
\sum_{h\in\left[  p\right]  }\left(  \sum_{i\in Z_{h},\text{ }j\in\left[
k\right]  \backslash Z_{h}}e_{B}\left(  V_{i},V_{j}\right)  \right)   &
=\sum_{h\in\left[  p\right]  }\left(  \sum_{s\in\left[  p\right]
\backslash\left\{  h\right\}  }\left(  \sum_{i\in Z_{h},\text{ }j\in Z_{s}%
}e_{B}\left(  V_{i},V_{j}\right)  \right)  \right) \\
&  =2\sum_{1\leq h<s\leq p}\left(  \sum_{i\in Z_{h},\text{ }j\in Z_{s}}%
e_{B}\left(  V_{i},V_{j}\right)  \right)  \geq\alpha N^{2}.
\end{align*}
Hence, we can select $h\in\left[  p\right]  $ so that
\[
\sum_{i\in Z_{h}}\left\{  e_{B}\left(  V_{i},V_{j}\right)  :j\in\left[
k\right]  \backslash Z_{h}\right\}  \geq\frac{\alpha N^{2}}{p}.
\]
Since by (\ref{maxd1}) we have%
\[
\left\vert Z_{h}\right\vert \leq\Delta\left(  B^{\ast}\right)  +1<\left(
1+2\xi\right)  \frac{k}{p}+1\leq\left(  1+3\xi\right)  \frac{k}{p},
\]
there is an $i_{1}\in Z_{h}$ such that
\[
\sum_{j\in\left[  k\right]  \backslash Z_{h}}e_{B}\left(  V_{i_{1}}%
,V_{j}\right)  \geq\frac{\alpha N^{2}}{\left(  1+3\xi\right)  k}\geq
\frac{\alpha N}{\left(  1+3\xi\right)  }t,
\]
and so,
\[
\sum_{j\in\left[  p\right]  \backslash\left\{  h\right\}  }\left(  \sum
_{\mu\left(  s\right)  =j}e_{B}\left(  V_{i_{1}},V_{s}\right)  \right)
\geq\frac{\alpha N}{\left(  1+3\xi\right)  }t
\]
Furthermore, in view of (\ref{minz}) we have%
\[
\sum_{j\in\left[  p\right]  \backslash\left\{  h\right\}  }\left\vert
Z_{j}\right\vert =k-\left\vert Z_{h}\right\vert \leq k-\left(  1-2p\xi\right)
\frac{k}{p}=\frac{p-1+2p\xi}{p}k,
\]
and thus%
\[
\frac{N}{\left(  1+3\xi\right)  }>\left(  \frac{p-1+2p\xi}{p}\right)
N\geq\left(  \frac{p-1+2p\xi}{p}\right)  kt\geq t\sum_{j\in\left[  p\right]
\backslash\left\{  h\right\}  }\left\vert Z_{j}\right\vert .
\]
Therefore,
\[
\sum_{j\in\left[  p\right]  \backslash\left\{  h\right\}  }\left(  \sum
_{\mu\left(  s\right)  =j}e_{B}\left(  V_{i_{1}},V_{s}\right)  \right)
\geq\alpha t^{2}\sum_{j\in\left[  p\right]  \backslash\left\{  h\right\}
}\left\vert Z_{j}\right\vert
\]
and the pigeonhole principle gives some $j\in\left[  p\right]  \backslash
\left\{  h\right\}  $ for which
\[
\sum_{\mu\left(  s\right)  =j}e_{B}\left(  V_{i_{1}},V_{s}\right)
>\alpha\left\vert Z_{j}\right\vert t^{2},
\]
completing the proof.
\end{proof}

\bigskip

\begin{proof}
[\textbf{Proof of Claim \ref{cl3.0}}]Set
\[
W_{0}=\left\{  u:u\in V_{i_{1}}\text{, }\left\vert \Gamma_{B}\left(  u\right)
\cap X\right\vert >\frac{\alpha}{2}\left\vert X\right\vert \right\}  .
\]
In view of of (\ref{cl1in}),
\begin{align*}
\alpha\left\vert X\right\vert t  &  <\sum_{u\in V_{i_{1}}}\left\vert
\Gamma_{B}\left(  u\right)  \cap X\right\vert =\sum_{u\in W_{0}}\left\vert
\Gamma_{B}\left(  u\right)  \cap X\right\vert +\sum_{u\in V_{i_{1}}\backslash
W_{0}}\left\vert \Gamma_{B}\left(  u\right)  \cap X\right\vert \\
&  <\left\vert W_{0}\right\vert \left\vert X\right\vert +\left(  t-\left\vert
W_{0}\right\vert \right)  \frac{\alpha}{2}\left\vert X\right\vert ,
\end{align*}
implying that
\[
\frac{\alpha}{2}t<\left(  1-\frac{\alpha}{2}\right)  \left\vert W_{0}%
\right\vert ,
\]
so $\left\vert W_{0}\right\vert >\left(  \alpha/2\right)  t.$
\end{proof}

\bigskip

\begin{proof}
[\textbf{Proof of Claim \ref{cl3.1}}]Let
\[
W=\left\{  u:u\in V_{i_{1}}\text{, }\left\vert \Gamma_{B}\left(  u\right)
\cap Y\right\vert >\left(  1-\sqrt{\beta}\right)  \left\vert Y\right\vert
\right\}
\]
We shall show that $\left\vert W\right\vert >\left(  1-\sqrt{\beta}\right)
t.$ Indeed,
\begin{align*}
\left(  1-\beta\right)  \left\vert Y\right\vert t  &  <\sum_{s\in
Z_{1}\backslash\left\{  i_{1}\right\}  }e\left(  V_{i_{1}},V_{s}\right)
=e\left(  V_{i_{1}},Y\right)  =\sum_{u\in V_{i_{1}}}\left\vert \Gamma
_{B}\left(  u\right)  \cap Y\right\vert \\
&  =\sum_{u\in W}\left\vert \Gamma_{B}\left(  u\right)  \cap Y\right\vert
+\sum_{u\in V_{i_{1}}\backslash W}\left\vert \Gamma_{B}\left(  u\right)  \cap
Y\right\vert \\
&  <\left\vert W\right\vert \left\vert Y\right\vert +\left(  t-\left\vert
W\right\vert \right)  \left(  1-\sqrt{\beta}\right)  \left\vert Y\right\vert .
\end{align*}
Hence
\[
\left(  1-\beta\right)  t<\left\vert W\right\vert +\left(  t-\left\vert
W\right\vert \right)  \left(  1-\sqrt{\beta}\right)  ,
\]
so $\left\vert W\right\vert >\left(  1-\sqrt{\beta}\right)  t.$

Now$W_{1}=W_{0}\cap W$ satisfies%
\[
\left\vert W_{1}\right\vert \geq\left\vert W_{0}\right\vert +\left\vert
W\right\vert -t>\left(  \frac{\alpha}{2}-\sqrt{\beta}\right)  t\geq
\frac{\alpha}{4}t,
\]
completing the proof.
\end{proof}

\bigskip

\begin{proof}
[\textbf{Proof of claim \ref{cl4}}]Let $\left\{  i_{1},\ldots,i_{s}\right\}  $
induces a maximal clique in $R^{\ast}$ containing $i_{1};$ assume for a
contradiction that $s<p.$ Then by (\ref{mindR}),
\[
d_{R^{\ast}}\left(  \left\{  i_{1},\ldots,i_{s}\right\}  \right)  \geq
\sum_{j=1}^{s}d_{R^{\ast}}\left(  i_{j}\right)  -\left(  s-1\right)
k>s\left(  \frac{p-1}{p}-2\xi\right)  k-\left(  s-1\right)  k>0.
\]
Thus, there is a vertex $i\in\left[  k\right]  $ joined in $R^{\ast}$ to all
vertices $i_{1},\ldots,i_{s},$ contradicting the fact that $\left\{
i_{1},\ldots,i_{s}\right\}  $ induces a maximal clique and completing the proof.
\end{proof}

\bigskip

\begin{proof}
[\textbf{Proof of Claim \ref{cl5}}]For $s=2,\ldots,p,$ applying Lemma
\ref{tKL}, select $P_{s}\subset V_{i_{1}}$ with $\left\vert P_{s}\right\vert
\geq\left(  1-\varepsilon\right)  t$ and $\left\vert \Gamma_{R}\left(
u\right)  \cap V_{i_{s}}\right\vert >\left(  \beta-\varepsilon\right)  t$ for
every $u\in P_{s}.$ Hence
\[
\left\vert
{\textstyle\bigcap_{s=2}^{p}}
P_{s}\right\vert >\left(  p-1\right)  \left(  1-\varepsilon\right)  t-\left(
p-2\right)  t>\left(  1-p\varepsilon\right)  t.
\]
Therefore, for $W_{2}=W_{1}\cap\left(
{\textstyle\bigcap_{s=2}^{p}}
P_{s}\right)  $ we have
\[
\left\vert W_{2}\right\vert =\left\vert W_{1}\cap\left(
{\textstyle\bigcap_{s=2}^{p}}
P_{s}\right)  \right\vert \geq\left\vert W_{1}\right\vert +\left\vert
{\textstyle\bigcap_{s=2}^{p}}
P_{s}\right\vert -t\geq\left(  \alpha/4\right)  t+\left(  1-p\varepsilon
\right)  t-t\geq\left(  \alpha/8\right)  t.
\]
Set $Q_{1}=W_{2}$ and let $a=a\left(  2,\beta/2,\gamma/\left(  2p\right)
\right)  $ (see Lemma \ref{probl}).

For $s=2,\ldots,p,$ applying Lemma \ref{probl} with $k=2,$ $d=\beta/2,$
$\lambda=\gamma/\left(  2p\right)  ,$ find $Q_{s}\subset Q_{s-1}$ with
$\left\vert Q_{s}\right\vert \geq\left\vert Q_{s-1}\right\vert ^{1-\gamma
/\left(  2p\right)  }$ and $\left\vert \Gamma_{R}\left(  uv\right)  \cap
V_{i_{s}}\right\vert >at$ for every $\left\{  u,v\right\}  \in Q_{s}^{\left(
2\right)  }.$ Set $W=Q_{p}$ and note that
\begin{align*}
\left\vert W\right\vert  &  =\left\vert Q_{p}\right\vert \geq\left\vert
Q_{p-1}\right\vert ^{1-\gamma/\left(  2p-2\right)  }\geq\cdots\geq\left\vert
Q_{1}\right\vert ^{\left(  1-\gamma/2p\right)  ^{p-1}}\geq\left\vert
Q_{1}\right\vert ^{1-\gamma\left(  p-1\right)  /\left(  2p\right)  }\\
&  \geq\left(  \frac{\alpha}{8}t\right)  ^{1-\gamma\left(  p-1\right)
/\left(  2p\right)  }>t^{1-\gamma/2},
\end{align*}
for $t$ sufficiently large.

Assume for a contradiction that $R\left[  W\right]  $ contains an edge $uv$.
Since $\left\vert \Gamma_{R}\left(  uv\right)  \cap V_{i_{s}}\right\vert >at,$
by Lemma \ref{ctL} we have
\begin{align*}
js_{p+1}\left(  R\right)   &  \geq\left(  \left(  a-\varepsilon\right)
t\right)  ^{p-1}>\left(  \left(  a-\varepsilon\right)  \frac{N\left(
1-\varepsilon\right)  }{K}\right)  ^{p-1}\\
&  >\left(  \frac{\left(  a-\varepsilon\right)  \left(  1-\varepsilon\right)
}{K}\right)  ^{p-1}N^{p-1}>cn^{p-1},
\end{align*}
a contradiction with (\ref{massum}). So $W$ is a clique in $B,$ completing the proof.
\end{proof}

\bigskip

\begin{proof}
[\textbf{Proof of Claim \ref{cl6}}]Since (\ref{minex}) and (\ref{miney}) imply
that $e_{B}\left(  X\right)  >\left(  1-\beta\right)  \left\vert X\right\vert
^{2}/2$ and $e_{B}\left(  Y\right)  >\left(  1-\beta\right)  \left\vert
Y\right\vert ^{2}/2,$ by Proposition \ref{dendeg}, there exist $X_{0}\subset
X$ and $Y_{0}\subset Y$ such that%
\begin{align*}
\left\vert X_{0}\right\vert  &  >\left(  1-\sqrt{\beta}\right)  \left\vert
X\right\vert >\left(  1-\sqrt{\beta}\right)  \left(  1-2p\xi\right)  \frac
{k}{p}t\geq\left(  1-3p\xi\right)  \frac{k}{p}t,\\
\delta\left(  B\left[  X_{0}\right]  \right)   &  >\left(  1-2\sqrt{\beta
}\right)  \left\vert X_{0}\right\vert ,\\
\left\vert Y_{0}\right\vert  &  >\left(  1-\sqrt{\beta}\right)  \left\vert
Y\right\vert >\left(  1-\sqrt{\beta}\right)  \left(  1-2p\xi\right)  \frac
{k}{p}t\geq\left(  1-3p\xi\right)  \frac{k}{p}t,\\
\delta\left(  B\left[  Y_{0}\right]  \right)   &  >\left(  1-2\sqrt{\beta
}\right)  \left\vert Y_{0}\right\vert .
\end{align*}
Also, for every $u\in W,$%
\begin{align*}
\left\vert \Gamma_{B}\left(  u\right)  \cap X_{0}\right\vert  &
\geq\left\vert \Gamma_{B}\left(  u\right)  \cap X\right\vert -\left\vert
X\backslash X_{0}\right\vert \geq\frac{\alpha}{4}\left\vert X\right\vert
-\sqrt{\beta}\left\vert X\right\vert >\frac{\alpha}{8}\left\vert
X_{0}\right\vert ,\\
\left\vert \Gamma_{B}\left(  u\right)  \cap Y_{0}\right\vert  &
\geq\left\vert \Gamma_{B}\left(  u\right)  \cap Y\right\vert -\left\vert
Y\backslash Y_{0}\right\vert \geq\left(  1-\sqrt{\beta}\right)  \left\vert
Y\right\vert -\sqrt{\beta}\left\vert Y\right\vert >\left(  1-2\sqrt{\beta
}\right)  \left\vert Y_{0}\right\vert .
\end{align*}

Next, Lemma \ref{probl} implies that there exists $a>0$ and $U\subset W$ such
that for every $Q\subset U^{\left(  2q\right)  },$ $\left\vert \Gamma
_{B}\left(  Q\right)  \cap X_{0}\right\vert >a\left\vert X_{0}\right\vert $
and $\left\vert U\right\vert >\left\vert W\right\vert ^{1-\gamma/2}.$

Also Lemma \ref{triml} implies that there exists $M\subset V\left(  H\right)
$ with $\left\vert M\right\vert \leq\left(  2q+1\right)  \left\vert
H\right\vert ^{1-\gamma}$ such that $\psi\left(  H-M\right)  \leq
\eta\left\vert H\right\vert $ and $d_{M}\left(  u\right)  \leq2q$ for every
$u\in V\left(  H-M\right)  .$ Since the graph $H\left[  M\right]  $ is
$q$-degenerate, for $t$ large, we have
\[
\left\vert M\right\vert \leq\left(  2q+1\right)  \left\vert H\right\vert
^{1-\gamma}<\left(  2q+1\right)  \left(  kt\right)  ^{1-\gamma}<t^{\left(
1-\gamma/2\right)  ^{2}}<\left\vert U\right\vert
\]
for $t$ large.

Let $\varphi:H\left[  M\right]  \rightarrow U$ be a one-to-one mapping; since
$B\left[  U\right]  $ is complete, $\varphi$ is a monomorphism. We shall
extend $\varphi$ to $H$ by mapping almost all components of $H-M$ into $Y_{0}$
and the remaining components into $X_{0}$. We can partition $H-M$ into two
disjoint graphs $H_{1}$ and $H_{2}$ such that
\begin{align}
\left\vert H_{1}\right\vert  &  <\left(  1-6q\sqrt{\beta}-3p\xi\right)
\frac{k}{p}t,\label{h1ub}\\
\left\vert H_{2}\right\vert  &  <\left(  a-2q\sqrt{\beta}-3p\xi\right)
\frac{k}{p}t. \label{h2ub}%
\end{align}
Indeed, collect into $H_{1}$ as many components of $H-M$ as possible so that
(\ref{h1ub}) still holds, and collect the remaining components into $H_{2}$.
Since $\psi\left(  H-M\right)  <\eta n,$ inequality (\ref{h2ub}) follows from%
\begin{align*}
\left\vert H_{2}\right\vert  &  \leq n-\left\vert H_{1}\right\vert \leq
n-\left(  1-6q\sqrt{\beta}-3p\xi\right)  \frac{k}{p}t+\eta n\\
&  <\left(  1+2\eta\right)  \frac{N}{p}-\left(  1-6q\sqrt{\beta}-3p\xi\right)
\frac{k}{p}t\\
&  <\left(  1+2\eta\right)  \frac{\left(  1+2\varepsilon\right)  k}%
{p}t-\left(  1-6q\sqrt{\beta}-3p\xi\right)  \frac{k}{p}t\\
&  <\left(  3\eta+2\varepsilon+6q\sqrt{\beta}+3p\xi\right)  \frac{k}%
{p}t<\left(  a-2q\sqrt{\beta}-3p\xi\right)  \frac{k}{p}t.
\end{align*}

Set $l=\left\vert H_{1}\right\vert ;$ Proposition \ref{arang} implies that the
vertices of $H_{1}$ can be arranged as $v_{1},\ldots,v_{l}$ so that
$\left\vert \Gamma_{H}\left(  v_{i}\right)  \cap\left\{  v_{1},\ldots
,v_{i-1}\right\}  \right\vert \leq q$ for every $i\in\left[  l\right]  .$ We
shall extend $\varphi$ over $H_{1}$ by mapping each $v_{i}\in V\left(
H_{1}\right)  $ in turn. Let $\Gamma_{H}\left(  v_{i}\right)  =\left\{
v_{i_{1}},\ldots,v_{i_{h}}\right\}  \cup\left\{  u_{i_{1}},\ldots,u_{i_{s}%
}\right\}  ,$ where $v_{i_{1}},\ldots,v_{i_{h}}\in\left\{  v_{1}%
,\ldots,v_{i-1}\right\}  ,$ $h\leq q,$ and $u_{i_{1}},\ldots,u_{i_{s}}\in M,$
$s\leq2q.$ Therefore,
\[
v_{i}\in\left(
{\textstyle\bigcap_{j=1}^{h}}
\Gamma_{H}\left(  v_{i_{j}}\right)  \right)  \cap\left(
{\textstyle\bigcap_{j=1}^{s}}
\Gamma_{H}\left(  u_{i_{j}}\right)  \right)  .
\]
Set for convenience $x_{j}=\varphi\left(  v_{i_{j}}\right)  $ for all
$j\in\left[  h\right]  ,$ and $y_{j}=\varphi\left(  u_{i_{j}}\right)  $ for
all $j\in\left[  s\right]  .$ Note that
\begin{align*}
&  \left(
{\textstyle\bigcap_{j=1}^{h}}
\Gamma_{B}\left(  x_{j}\right)  \cap Y_{0}\right)  \cap\left(
{\textstyle\bigcap_{j=1}^{s}}
\Gamma_{B}\left(  y_{j}\right)  \cap Y_{0}\right) \\
&  \geq\sum_{j\in\left[  h\right]  }\left\vert \Gamma_{B}\left(  x_{j}\right)
\cap Y_{0}\right\vert +\sum_{j\in\left[  s\right]  }\left\vert \Gamma
_{B}\left(  y_{j}\right)  \cap Y_{0}\right\vert -\left(  h+s-1\right)
\left\vert Y_{0}\right\vert \\
&  >\left(  h+s\right)  \left(  1-2\sqrt{\beta}\right)  \left\vert
Y_{0}\right\vert -\left(  h+s-1\right)  \left\vert Y_{0}\right\vert >\left(
1-6q\sqrt{\beta}\right)  \left\vert Y_{0}\right\vert \\
&  >\left(  1-6q\sqrt{\beta}\right)  \left(  1-3p\xi\right)  \frac{k}%
{p}>\left\vert H_{1}\right\vert .
\end{align*}
Hence, there is a vertex $z\in Y_{0}$ that is joined to the vertices
$x_{1},\ldots,x_{h},y_{1},\ldots,y_{s}$ and is outside the current range of
$\varphi.$ Setting $\varphi\left(  v_{i}\right)  =z,$ we extend $\varphi$ to a
monomorphism that maps $v_{i}$ into $Y_{0}$ as well. In this way $\varphi$ can
be extended over the entire $H_{1}$.

Set now $l=\left\vert H_{2}\right\vert ;$ Proposition \ref{arang} implies that
the vertices of $H_{2}$ can be arranged as $v_{1},\ldots,v_{l}$ so that
$\left\vert \Gamma_{H}\left(  v_{i}\right)  \cap\left\{  v_{1},\ldots
,v_{i-1}\right\}  \right\vert \leq q$ for every $i\in\left[  l\right]  .$ We
shall extend $\varphi$ over $H_{2}$ mapping each $v_{i}\in V\left(
H_{2}\right)  $ in turn. Let $\Gamma_{H}\left(  v_{i}\right)  =\left\{
v_{i_{1}},\ldots,v_{i_{h}}\right\}  \cup\left\{  u_{i_{1}},\ldots,u_{i_{s}%
}\right\}  $ where $v_{i_{1}},\ldots,v_{i_{s}}\in\left\{  v_{1},\ldots
,v_{i-1}\right\}  ,$ $h\leq q,$ and $u_{i_{1}},\ldots,u_{i_{s}}\in M,$
$s\leq2q.$ Therefore,
\[
v_{i}\in\left(
{\textstyle\bigcap_{j=1}^{h}}
\Gamma_{H}\left(  v_{i_{j}}\right)  \right)  \cap\left(
{\textstyle\bigcap_{j=1}^{s}}
\Gamma_{H}\left(  u_{i_{j}}\right)  \right)  .
\]
Set for convenience $x_{j}=\varphi\left(  v_{i_{j}}\right)  $ for all
$j\in\left[  h\right]  ,$ and $y_{j}=\varphi\left(  u_{i_{j}}\right)  $ for
all $j\in\left[  s\right]  .$ Note that%
\begin{align*}
&  \left(
{\textstyle\bigcap_{j=1}^{h}}
\Gamma_{B}\left(  x_{j}\right)  \cap X_{0}\right)  \cap\left(
{\textstyle\bigcap_{j=1}^{s}}
\Gamma_{B}\left(  y_{j}\right)  \cap X_{0}\right) \\
&  \geq a\left\vert X_{0}\right\vert +\sum_{j\in\left[  s\right]  }\left\vert
\Gamma_{B}\left(  y_{j}\right)  \cap X_{0}\right\vert -s\left\vert
X_{0}\right\vert \\
&  >a\left\vert X_{0}\right\vert +s\left(  1-2\sqrt{\beta}\right)  \left\vert
X_{0}\right\vert -s\left\vert X_{0}\right\vert \\
&  >\left(  a-2q\sqrt{\beta}\right)  \left\vert X_{0}\right\vert >\left(
a-2q\sqrt{\beta}\right)  \left(  1-3p\xi\right)  \frac{k}{p}>\left\vert
H_{2}\right\vert .
\end{align*}
Hence, there is a vertex $z\in X_{0}$ that is joined to the vertices
$x_{1},\ldots,x_{h},y_{1},\ldots,y_{s}$ and is outside the current range of
$\varphi.$ Setting $\varphi\left(  v_{i}\right)  =z,$ we extend $\varphi$ to a
monomorphism that maps $v_{i}$ into $X_{0}$ as well. In this way $\varphi$ can
be extended over the entire $H_{2}$.
\end{proof}

\bigskip

\begin{proof}
[\textbf{Proof of Claim \ref{cl7}}]In view of (\ref{loblu}) and (\ref{minz}),
\begin{align*}
e\left(  R\right)   &  \geq\sum_{1\leq h<s\leq p}\left(  \sum_{i\in
Z_{h},\text{ }j\in Z_{s}}e_{R}\left(  V_{i},V_{j}\right)  \right)  \geq
\sum_{1\leq h<s\leq p}\left\vert Z_{h}\right\vert \left\vert Z_{s}\right\vert
t^{2}-\sum_{1\leq h<s\leq p}\left(  \sum_{i\in Z_{h},\text{ }j\in Z_{s}}%
e_{B}\left(  V_{i},V_{j}\right)  \right) \\
&  \geq\binom{p}{2}\left(  1-2p\xi\right)  ^{2}\frac{k^{2}t^{2}}{p^{2}}%
-\frac{\alpha}{2}N^{2}=\frac{p-1}{2p}\left(  1-4p\xi\right)  \left(
1-\varepsilon\right)  ^{2}N^{2}-\frac{\alpha}{2}N^{2}\\
&  \geq\left(  \frac{p-1}{2p}-4p\xi-2\varepsilon-\frac{\alpha}{2}\right)
N^{2}\geq\left(  \frac{p-1}{2p}-\alpha\right)  N^{2}.
\end{align*}
On the other hand we have
\[
js_{p+1}\left(  R\right)  <cn^{p-1}<\left(  1-\frac{1}{p^{3}}\right)
\frac{N^{p-1}}{p^{p+5}};
\]
hence, Fact \ref{stabj} implies that $R$ has a $p$-partite induced subgraph
$R_{0}$ with $\left\vert R_{0}\right\vert >\left(  1-2\sqrt{\alpha}\right)  N$
and%
\begin{equation}
\delta\left(  R_{0}\right)  >\left(  1-\frac{1}{p}-4\sqrt{\alpha}\right)  N.
\label{mindg1}%
\end{equation}

We shall find $R_{1}$ as an induced subgraph of $R_{0}.$ Observe that by
(\ref{mindg1}) every color class of $R_{0}$ has at most $N-\delta\left(
R_{0}\right)  >\left(  1/p+4\sqrt{\alpha}\right)  N$ vertices. Hence, every
color class of $G_{0}$ has at least
\[
\left(  1-2\sqrt{\alpha}\right)  N-\left(  p-1\right)  \left(  \frac{1}%
{p}+4\sqrt{\alpha}\right)  N>\left(  1-4p\left(  p-1\right)  \sqrt{\alpha
}\right)  \frac{N}{p}>\left(  1-\theta\right)  n
\]
vertices. From each color class select a set of $\left\lceil \left(
1-\theta\right)  n\right\rceil $ vertices and write $R_{1}$ for the graph
induced by their union.

Let $u\in V\left(  R_{1}\right)  $ and $U$ be a color class of $R_{1}$ such
that $u\notin U.$ Since $\delta\left(  R_{1}\right)  \geq\delta\left(
R_{0}\right)  -\left\vert R_{0}\right\vert +\left\vert R_{1}\right\vert ,$ we
see that
\begin{align*}
\left\vert \Gamma_{R_{1}}\left(  u\right)  \cap U\right\vert  &  >\left\vert
U\right\vert +\delta\left(  R_{1}\right)  -\frac{p-1}{p}\left\vert
R_{1}\right\vert \geq\left\vert U\right\vert +\delta\left(  R_{0}\right)
-\left\vert R_{0}\right\vert +\left\vert R_{1}\right\vert -\frac{p-1}%
{p}\left\vert R_{1}\right\vert \\
&  =\delta\left(  R_{0}\right)  -\left\vert R_{0}\right\vert +\frac{2}%
{p}\left\vert R_{1}\right\vert >\left(  1-\frac{1}{p}-4\sqrt{\alpha}\right)
N-N+\left(  \frac{2}{p}-8p\sqrt{\alpha}\right)  N\\
&  >\left(  1-8p\left(  p+1\right)  \sqrt{\alpha}\right)  \frac{N}{p}%
\geq\left(  1-\theta\right)  n,
\end{align*}
completing the proof.
\end{proof}

\bigskip

\begin{proof}
[\textbf{Proof of Claim \ref{cl8}}]Set $s=\left\vert U_{1}\right\vert
=\cdots=\left\vert U_{p}\right\vert .$ According to Claim \ref{cl7},
\begin{align}
\left(  1-\theta\right)  n  &  <s<n,\label{bnds}\\
\left\vert \Gamma_{R}\left(  u\right)  \cap U_{i}\right\vert  &  >\left(
1-\theta\right)  n\nonumber
\end{align}
for every $U_{i}$ and every $u\in V\left(  R_{1}\right)  \backslash U_{i}$.

Set $X=V\left(  R\right)  \backslash V\left(  R_{1}\right)  $ and define a
partition $X=Y\cup Z$ as follows:
\begin{align*}
Y  &  =\left\{  u:u\in X,\text{ }\Gamma_{R}\left(  u\right)  \cap U_{i}%
\neq\varnothing\text{ for every }i\in\left[  p\right]  \right\}  ,\\
Z  &  =X\backslash Y.
\end{align*}

We first show that for every $u\in Y,$ there exists two distinct color classes
$U_{i}$ and $U_{j}$ such that
\begin{equation}
\left\vert \Gamma_{R}\left(  u\right)  \cap U_{i}\right\vert \leq p^{2}\theta
n,\text{ \ \ \ }\left\vert \Gamma_{R}\left(  u\right)  \cap U_{j}\right\vert
\leq p^{2}\theta n. \label{mingR}%
\end{equation}
For a contradiction, assume the opposite: let $u\in Y$ be such that
$\left\vert \Gamma_{R}\left(  u\right)  \cap U_{i}\right\vert >\theta n$ for
at least $p-1$ values $i\in\left[  p\right]  ,$ say for $i=2,\ldots,p.$ The
definition of $Y$ implies that there exists some $v\in U_{1}\cap\Gamma
_{R}\left(  u\right)  .$ Observe that for every $i\in\left[  2..p\right]  ,$%
\begin{align*}
\left\vert \Gamma_{R}\left(  u\right)  \cap\Gamma_{R}\left(  v\right)  \cap
U_{i}\right\vert  &  \geq\left\vert \Gamma_{R}\left(  u\right)  \cap
U_{i}\right\vert +\left\vert \Gamma_{R}\left(  v\right)  \cap U_{i}\right\vert
-\left\vert U_{i}\right\vert \\
&  >p^{2}\theta n+n-\theta n-s>\left(  p^{2}-1\right)  \theta n.
\end{align*}
Therefore, for every $i\in\left[  2..p\right]  ,$ we can select a set
$W_{i}\subset\Gamma_{R}\left(  u\right)  \cap\Gamma_{R}\left(  v\right)  \cap
U_{i}$ with
\[
\left\vert W_{i}\right\vert =m=\left\lceil \left(  p^{2}-1\right)  \theta
n\right\rceil .
\]
We shall prove that the set $W=\cup_{i=2}^{p}W_{i}$ induces at least
\[
\frac{1}{\left(  p-1\right)  ^{2}}\left(  \left(  p^{2}-1\right)
\theta\right)  ^{p-1}n^{p-1}%
\]
$\left(  p-1\right)  $-cliques in $R$ and thus obtain a contradiction with
(\ref{massum}). The assertion is immediate for $p=2;$ assume henceforth that
$p\geq3.$ Let $w\in W$ be a vertex of minimum degree in $R\left[  W\right]  ,$
say let $w\in W_{i}.$ We have
\begin{align*}
\delta\left(  R\left[  W\right]  \right)   &  =\sum_{j\in\left[  2..p\right]
\backslash\left\{  i\right\}  }\left\vert \Gamma_{R}\left(  w\right)  \cap
W_{j}\right\vert \geq\sum_{j\in\left[  2..p\right]  \backslash\left\{
i\right\}  }\left\vert \Gamma_{R}\left(  w\right)  \cap U_{j}\right\vert
+\left\vert W_{j}\right\vert -\left\vert U_{j}\right\vert \\
&  >\left(  p-2\right)  \left(  \left(  1-\theta\right)  n+m-n\right)
=\left(  p-2\right)  \left(  m-\theta n\right) \\
&  \geq\left(  p-2\right)  \left(  1-\frac{1}{p^{2}-1}\right)  m.
\end{align*}

Hence, in view of $\left\vert W\right\vert =\left(  p-1\right)  m,$
\[
\delta\left(  R\left[  W\right]  \right)  >\frac{p-2}{p-1}\left(  1-\frac
{1}{p^{2}-1}\right)  \left\vert W\right\vert =\left(  \frac{p-3}{p-2}%
+\frac{4p-5}{\left(  p-1\right)  ^{2}\left(  p+1\right)  \left(  p-2\right)
}\right)  \left\vert W\right\vert .
\]

Applying Fact \ref{minj} to $R\left[  W\right]  ,$ we obtain
\begin{align*}
js_{p+1}\left(  R\right)   &  \geq k_{p-1}\left(  R\left[  W\right]  \right)
>\frac{4p-5}{\left(  p-1\right)  ^{2}\left(  p+1\right)  \left(  p-2\right)
}\cdot\frac{\left(  p-2\right)  ^{2}}{\left(  p-1\right)  }\left(
\frac{\left\vert W\right\vert }{p-1}\right)  ^{p-1}\\
&  \geq\frac{1}{\left(  p-1\right)  ^{2}}\left(  \frac{\left\vert W\right\vert
}{p-1}\right)  ^{p-1}\geq\frac{1}{\left(  p-1\right)  ^{2}}\left(  \left(
p^{2}-1\right)  \theta\right)  ^{p-1}n^{p-1}>cn^{p-2},
\end{align*}
a contradiction with (\ref{massum}).

Hence, for every $u\in Y,$ there exists two distinct color classes $U_{i}$ and
$U_{j}$ such that (\ref{mingR}) holds. For every $i\in\left[  p\right]  ,$
set
\begin{align}
Z_{i}  &  =\left\{  u:u\in Z,\text{ }\Gamma_{R}\left(  u\right)  \cap
U_{i}=\varnothing\right\} \label{prop2i}\\
Y_{i}  &  =\left\{  u:u\in Y,\text{ }\Gamma_{R}\left(  u\right)  \cap
U_{i}\leq p^{2}\theta\right\}  \label{prop3i}%
\end{align}
We have
\[
\sum_{i=1}^{p}\left\vert Z_{i}\right\vert \geq\left\vert \cup_{i=1}^{p}%
Z_{i}\right\vert =\left\vert Z\right\vert ,\text{ \ \ and \ \ }\sum_{i=1}%
^{p}\left\vert Y_{i}\right\vert \geq2\left\vert \cup_{i=1}^{p}Y_{i}\right\vert
=2\left\vert Y\right\vert .
\]
Hence
\[
\sum_{i=1}^{p}\left(  \left\vert U_{i}\right\vert +\left\vert Z_{i}\right\vert
+\left\vert Y_{i}\right\vert \right)  \geq N+\left\vert Y\right\vert =p\left(
n-1\right)  +1+\left\vert Y\right\vert ,
\]
and there exists $i\in\left[  p\right]  $ such that
\[
\left\vert U_{i}\right\vert +\left\vert Z_{i}\right\vert +\left\vert
Y_{i}\right\vert \geq n-1+\left\lceil \frac{\left\vert Y_{i}\right\vert +1}%
{p}\right\rceil .
\]

Set $M=U_{i},$ $A=Z_{i},$ $C=Y_{i}$ and apply the following procedure to the
sets $A$ and $C.$

\emph{While} $\left\vert M\right\vert +\left\vert A\right\vert +\left\vert
C\right\vert >n-1+\left\lceil \left(  \left\vert C\right\vert +1\right)
/p\right\rceil $ \emph{do }

\qquad\emph{if} $C\neq\varnothing$ \emph{remove a vertex from} $C$ \emph{else
remove a vertex from} $A;$

\emph{end.}

This procedure is defined correctly in view of $\left\vert M\right\vert =s$
and inequalities (\ref{bnds}). Upon the end of the procedure we have
\[
\left\vert M\right\vert +\left\vert A\right\vert +\left\vert C\right\vert
=n-1+\left\lceil \frac{\left\vert C\right\vert +1}{p}\right\rceil ,
\]
so condition (\ref{cond}) holds. We also see that
\[
\left\vert A\right\vert =n-1+\left\lceil \frac{\left\vert C\right\vert +1}%
{p}\right\rceil -\left\vert M\right\vert -\left\vert C\right\vert \leq
n-\left\vert M\right\vert <\theta n,
\]
so condition (\ref{cond2}) holds as well. Finally, condition (\ref{cond3})
holds due to
\begin{align*}
\frac{1}{2}\left\vert C\right\vert  &  \leq\frac{p-1}{p}\left\vert
C\right\vert =\left\vert C\right\vert -\frac{\left\vert C\right\vert +1}%
{p}-\frac{p-1}{p}+1\leq\left\vert C\right\vert -\left\lceil \frac{\left\vert
C\right\vert +1}{p}\right\rceil +1\\
&  \leq n-\left\vert M\right\vert <\theta n.
\end{align*}

To complete the proof of the claim, observe that property \emph{(i) }holds
since the set $M$ is independent in $R;$ properties \emph{(ii) }and
\emph{(iii) }hold in view of (\ref{prop2i}) and (\ref{prop3i}).
\end{proof}

\bigskip

\begin{proof}
[\textbf{Proof of Claim \ref{cl9}}]Define a set $M^{\prime}\subset M$ by
\[
M^{\prime}=\left\{  u:u\in M,\text{ }\left\vert \Gamma_{R}\left(  u\right)
\cap C\right\vert \geq\left(  1-2p^{2}\theta\right)  \left\vert C\right\vert
\right\}  ;
\]
first we shall prove that $\left\vert M^{\prime}\right\vert \geq\left\vert
M\right\vert /2.$ This is obvious if $C=\varnothing,$ so we shall assume that
$\left\vert C\right\vert >0.$ We have
\begin{align*}
\left(  1-p^{2}\theta\right)  \left\vert C\right\vert \left\vert M\right\vert
&  \leq\sum_{u\in C}\left\vert \Gamma_{B}\left(  u\right)  \cap M\right\vert
=e_{B}\left(  M,C\right)  =\sum_{u\in M}\left\vert \Gamma_{B}\left(  u\right)
\cap C\right\vert \\
&  =\sum_{u\in M^{\prime}}\left\vert \Gamma_{B}\left(  u\right)  \cap
C\right\vert +\sum_{u\in M\backslash M^{\prime}}\left\vert \Gamma_{B}\left(
u\right)  \cap C\right\vert \\
&  \leq\left\vert C\right\vert \left\vert M^{\prime}\right\vert +\left(
1-2p^{2}\theta\right)  \left\vert C\right\vert \left(  \left\vert M\right\vert
-\left\vert M^{\prime}\right\vert \right)  ,
\end{align*}
implying that
\[
\left(  1-p^{2}\theta\right)  \left\vert M\right\vert \leq\left\vert
M^{\prime}\right\vert +\left(  1-2p^{2}\theta\right)  \left(  \left\vert
M\right\vert -\left\vert M^{\prime}\right\vert \right)  =\left(
1-2p^{2}\theta\right)  \left\vert M\right\vert +2p^{2}\theta\left\vert
M^{\prime}\right\vert ,
\]
and the desired inequality follows.

Setting
\[
W_{0}=\left\{  u:u\in V\left(  H\right)  ,\text{ }d\left(  u\right)
\leq2q\right\}  ,
\]
by Proposition \ref{cheb} we have $\left\vert W_{0}\right\vert \geq n/\left(
2q+1\right)  .$ Since by Fact \ref{color} $H\left[  W_{0}\right]  $ is
$\left(  q+1\right)  $-partite, there exists an independent set $W_{1}\subset
W_{0}$ with
\[
\left\vert W_{1}\right\vert \geq\frac{\left\vert W_{0}\right\vert }{q+1}%
\geq\frac{n}{\left(  q+1\right)  \left(  2q+1\right)  }>\theta n.
\]

If $\left\vert A\right\vert +\left\vert M\right\vert \geq n,$ we map $H$ into
$M\cup A$ as follows:

- select a set $W\subset W_{1}$ with $\left\vert W\right\vert =\left\vert
A\right\vert $ - this is possible since $\left\vert A\right\vert <\theta n;$

- map arbitrarily $W$ into $A;$

- map arbitrarily $V\left(  H\right)  \backslash W$ into $M.$

This mapping is a monomorphism since the set $W$ is independent in $H,$ the
set $M$ induces a complete graph in $B$, and the sets $A$ and $M$ induce a
complete bipartite graph in $B.$

We assume henceforth that $\left\vert M\right\vert +\left\vert A\right\vert
<n.$ Since
\[
\left\vert C\right\vert <2\theta n\leq\frac{n}{\left(  q+1\right)  \left(
2q+1\right)  },
\]
select an independent set $W\subset W_{1}$ with
\[
\left\vert W\right\vert =n-\left\vert A\right\vert -\left\vert M\right\vert
=\left\vert C\right\vert +1-\left\lceil \frac{\left\vert C\right\vert +1}%
{p}\right\rceil ,
\]
and set $P=\cup_{u\in W}\Gamma_{H}\left(  u\right)  .$ Clearly
\[
\left\vert P\right\vert \leq\sum_{u\in W}\left\vert \Gamma_{H}\left(
u\right)  \right\vert \leq2q\left\vert W\right\vert \leq2q\left\vert
C\right\vert n\leq\frac{4q\theta}{1-\theta}\left\vert M\right\vert \leq
\frac{\left\vert M\right\vert }{2}\leq\left\vert M^{\prime}\right\vert .
\]

We construct a monomorphism $\varphi:$ $H\rightarrow B$ in two steps:
\emph{(a)} define $\varphi$ on $H\left[  W\cup P\right]  ;$ \emph{(b)} extend
$\varphi$ over $H-W-P.$

\emph{(a) }\textbf{defining a monomorphism }$\varphi:H\left[  W\cup P\right]
\rightarrow B$

Define $\varphi$ as an arbitrary one-to-one mapping $\varphi:P\rightarrow
M^{\prime}$ and extend $\varphi$ by mapping $W$ into $C$ one vertex at a time.
Suppose $W^{\prime}\subset W$ is the set of vertices already mapped; if
$W^{\prime}\neq W,$ select an unmapped $u\in W$ and let
\[
\left\{  v_{1},\ldots,v_{r}\right\}  =\varphi\left(  \Gamma_{H}\left(
u\right)  \right)  \subset M^{\prime}.
\]
Since $W\subset W_{1}\subset W_{0},$ we see that $r\leq2q.$ Then
\begin{align*}
\left\vert
{\textstyle\bigcap_{i=1}^{r}}
\left(  \Gamma_{B}\left(  v_{i}\right)  \cap C\right)  \right\vert  &
\geq\sum_{i=1}^{r}\left\vert \Gamma_{B}\left(  v_{i}\right)  \cap C\right\vert
-\left(  r-1\right)  \left\vert C\right\vert \\
&  \geq r\left\lceil \left(  1-2p^{2}\theta\right)  \left\vert C\right\vert
\right\rceil -\left(  r-1\right)  \left\vert C\right\vert \\
&  \geq\left\vert C\right\vert +r\left\lceil -2p^{2}\theta\left\vert
C\right\vert \right\rceil \geq\left\vert C\right\vert +\left\lceil
-4p^{2}q\theta\left\vert C\right\vert \right\rceil \geq\left\vert C\right\vert
+\left\lceil -\frac{\left\vert C\right\vert }{p}\right\rceil \\
&  =\left\vert C\right\vert +1-\left\lceil \frac{\left\vert C\right\vert
+1}{p}\right\rceil =\left\vert W\right\vert >\left\vert W^{\prime}\right\vert
.
\end{align*}
Hence, there exists a vertex $v\in\left(
{\textstyle\bigcap_{i=1}^{r}}
\left(  \Gamma_{B}\left(  v_{i}\right)  \cap C\right)  \right)  \backslash
\varphi\left(  W^{\prime}\right)  .$ Letting $v=\varphi\left(  u\right)  ,$ we
extend $\varphi$ to $W^{\prime}\cup\left\{  u\right\}  ;$ this extension can
be continued the entire $W$ is mapped into $C.$

\emph{(b) }\textbf{extending }$\varphi$ \textbf{over }$H-W-P$

Since $H-W-P$ is $\left(  q+1\right)  $-partite, the set $V\left(  H\right)
\backslash\left(  W\cup P\right)  $ contains an independent set $W^{\prime
\prime}$ with
\[
\left\vert W^{\prime\prime}\right\vert \geq\frac{n-\left\vert W\right\vert
-\left\vert P\right\vert }{q+1}\geq\frac{n-\left(  2q+1\right)  \left\vert
W\right\vert }{q+1}\geq\frac{1-\left(  2q+1\right)  \theta}{q+1}n\geq\theta
n>\left\vert A\right\vert .
\]

Now, extend $\varphi$ to $H$ by mapping arbitrarily $W^{\prime\prime}$ into
$A$ and $V\left(  H\right)  \backslash\left(  W\cup P\cup W^{\prime\prime
}\right)  $ into $M\backslash\varphi\left(  P\right)  .$ This extension is a
monomorphism due to the following facts:

- $W^{\prime\prime}$ is independent in $H,$

- the set $E_{H}\left(  W,W^{\prime\prime}\right)  $ is empty,

- the set $M$ induces a complete graph in $B,$

- the sets $A$ and $M$ induce a complete bipartite graph in $B.$

This completes the proof of the claim.
\end{proof}

\subsection{\label{pmth1}Proof of Theorem \ref{mth1}}

The proof of Theorem \ref{mth1} is reduced to the following proposition.

\begin{proposition}
\label{cor1}For every $p\geq3,$ $c>0,$ there exists $b>0$ such that if
$G=G\left(  n\right)  $ is a graph with $js_{p}\left(  G\right)  >cn^{p-2}$,
then $K_{p}\left(  1,1,t,\ldots,t\right)  \subset G,$ for $t>b\log n.$%
\hfill{$\Box$}
\end{proposition}

In turn, Proposition \ref{cor1} is implied by the following fact.

\begin{fact}
For every $p\geq3,$ $c>0$ there exists $b>0$ such that if $G=G\left(
n\right)  $ is a graph with $k_{p}\left(  G\right)  \geq cn^{p}$, then
$K_{p}\left(  t\right)  \subset G$ for $t\geq b\log n.$\hfill{$\Box$}
\end{fact}

The proof of this theorem can be found in \cite{Nik07}.

\subsection{\label{pmth2}Proof of Theorem \ref{mth2}}

\begin{lemma}
\label{le3}For every $p\geq2,$ $d\geq1$ and $c>0,$ there exists $\alpha>0,$
such that if $G=G\left(  n\right)  $ and $k_{p}\left(  G\right)  >cn^{p}$,
then $G$ contains every $p$-partite graph $H$ with $\left\vert H\right\vert
\leq\alpha n$ and $\Delta\left(  H\right)  \leq d.$
\end{lemma}

\begin{proof}
We sketch a proof using the Blow-up Lemma, see \cite{KSS97}. Applying the
Regularity Lemma of Szemer\'{e}di we first find an $\varepsilon$-regular
partition $V\left(  G\right)  =\cup_{i=0}^{k}V_{i}$ with $\varepsilon
\ll\left(  p,c\right)  $, $1/\varepsilon\leq k\leq K\left(  \varepsilon
\right)  .$ Remove the vertices from $V_{0}$ and all edges that belong to:

- any $E\left(  V_{i}\right)  ;$

- any irregular pair $\left(  V_{i},V_{j}\right)  ;$

- any pair $\left(  V_{i},V_{j}\right)  $ with $\sigma_{G}\left(  V_{i}%
,V_{j}\right)  <c.$

A straightforward counting shows that the remaining graph contains a $K_{p},$
and so there exists $p$ sets $V_{i_{1}},\ldots,V_{i_{p}}$ such that, for every
$1\leq l<j\leq p,$ the pair $\left(  V_{i_{l}},V_{i_{j}}\right)  $ is
$\varepsilon$-regular and $\sigma\left(  V_{i_{l}},V_{i_{j}}\right)  >c$.
Using Fact \ref{tKL}, we find subsets $U_{i_{j}}\subset V_{i_{j}}$ such that

- $\left\vert U_{i_{1}}\right\vert =\cdots=\left\vert U_{i_{p}}\right\vert
\geq\left(  1-p\varepsilon\right)  \left\vert V_{i_{1}}\right\vert ,$

- for every $1\leq l<j\leq p,$ the pair $\left(  U_{i_{l}},U_{i_{j}}\right)  $
is $2\varepsilon$-regular and every vertex $u\in U_{i_{l}}$ has at least $c/2$
neighbors in $U_{i_{j}}.$

According to the Blow-up Lemma, the graph $G\left[  \cup_{j=1}^{p}U_{i_{j}%
}\right]  $ contains all spanning graphs with maximum degree at most $d,$ for
$\left\vert U_{i_{1}}\right\vert $ sufficiently large. Therefore, $G\left[
\cup_{j=1}^{p}U_{i_{j}}\right]  $ contains all $p$-partite graphs of order
$\left\vert U_{i_{1}}\right\vert +p-1$ and of maximum degree at most $d.$
Since $\left\vert U_{i_{1}}\right\vert >n/\left(  2K\right)  ,$ the assertion follows.
\end{proof}

\subsection{\label{pprobl}Probabilistic Lemmas}

We deduce Lemma \ref{probl} from a more general result; its proof is an
adaptation of Sudakov's proof of Lemma 2.1 in \cite{Sud05}.

\begin{lemma}
\label{probl1}Suppose $G$ is a bipartite graph with parts $V$ and $U$ with
$\left\vert V\right\vert =n,$ $\left\vert U\right\vert =m,$ and $e\left(
G\right)  \geq dnm$. Let $H$ be a uniform $k$-graph with $V\left(  H\right)
=V$ and $d\left(  v_{1},\ldots,v_{k}\right)  \leq am$ for every $\left\{
v_{1},\ldots,v_{k}\right\}  \in E\left(  H\right)  .$ Then there exists
$W\subset V$ with $\left\vert W\right\vert \geq\left(  d^{i}/2\right)  n$ such
that $e\left(  H\left[  W\right]  \right)  \leq\left(  a/d\right)  ^{i}%
n^{k-1}\left\vert W\right\vert .$
\end{lemma}

\begin{proof}
Chose $I\in U^{i}$ uniformly. Let $W=\Gamma\left(  I\right)  $ and define the
random variables
\[
X=\left\vert W\right\vert ,\text{ \ \ }Y=e\left(  H\left[  W\right]  \right)
,\text{ \ \ }Z=X-\frac{d^{i}}{a^{i}n^{k-1}}Y-\frac{d^{i}}{2}n.
\]
We have
\begin{align*}
\mathbb{E}\left(  X\right)   &  =\frac{1}{m^{i}}\sum_{v\in V}d^{i}\left(
v\right)  \geq\frac{n}{m^{i}}\left(  \sum_{v\in V}\frac{d\left(  v\right)
}{n}\right)  ^{i}\geq\frac{n}{m^{i}}\left(  dm\right)  ^{i}=d^{i}n,\\
\mathbb{E}\left(  Y\right)   &  \leq\frac{1}{m^{i}}\sum_{\left\{  v_{1}%
,\ldots,v_{k}\right\}  \in E\left(  H\right)  }d^{i}\left(  v_{1},\ldots
,v_{k}\right)  \leq\frac{1}{m^{i}}e\left(  H\right)  \left(  am\right)
^{i}\leq a^{i}\frac{n^{k}}{2}\\
\mathbb{E}\left(  Z\right)   &  =\mathbb{E}\left(  X\right)  -\frac{d^{i}%
}{a^{i}n^{k-1}}\mathbb{E}\left(  Y\right)  -\frac{d^{i}}{2}n\geq\frac{d^{i}%
}{2}n-\frac{d^{i}}{a^{i}n^{k-1}}a^{i}\frac{n^{k}}{2}=0
\end{align*}
Thus, there exists $I_{0}\in U^{i}$ for which $\mathbb{E}\left(  Z\right)
\geq0.$ Then for $W=\Gamma\left(  I_{0}\right)  $ we have
\begin{align*}
\left\vert W\right\vert -\frac{d^{i}}{2}n  &  =X-\frac{d^{i}}{2}%
n=Z+\frac{d^{i}}{a^{i}n^{k-1}}Y\geq0,\\
e\left(  H\left[  W\right]  \right)   &  =Y=\frac{a^{i}n^{k-1}}{d^{i}}\left(
X-Z\right)  \leq\left(  \frac{a}{d}\right)  ^{i}n^{k-1}\left\vert W\right\vert
,
\end{align*}
completing the proof.
\end{proof}

\begin{proof}
[\textbf{Proof of Lemma \ref{probl}}]Set $a=d^{2k/\lambda+1}$ and
$n=\left\vert U_{1}\right\vert ;$ let $i$ be the smallest integer such that
$\left(  a/d\right)  ^{i}n^{k}<1,$ i.e.,
\[
i-1<\frac{k}{\ln\left(  d/a\right)  }\ln n=\frac{-\lambda}{2\ln d}\ln n.
\]
Define a $k$-uniform graph $H$ with $V\left(  H\right)  =U_{1}$: a $k$-set
$\left\{  u_{1},\ldots,u_{k}\right\}  \subset U_{1}$ belongs to $E\left(
H\right)  $ if $d\left(  u_{1},\ldots,u_{k}\right)  \leq a\left\vert
U_{2}\right\vert .$ According to Lemma \ref{probl1}, there exists $W\subset
U_{1}$ with $\left\vert W\right\vert \geq\left(  d^{i}/2\right)  n$ and
\[
e\left(  H\left[  W\right]  \right)  \leq\left(  \frac{a}{d}\right)
^{i}n^{k-1}\left\vert W\right\vert \leq\left(  \frac{a}{d}\right)  ^{i}%
n^{k}<1.
\]
Thus, $W$ is an independent set in $H,$ and so $d\left(  u_{1},\ldots
,u_{k}\right)  >a\left\vert U_{2}\right\vert $ for every $k$-set $\left\{
u_{1},\ldots,u_{k}\right\}  \subset W$. We also have, for $n$ large,
\[
\left\vert W\right\vert \geq\frac{d^{i}}{2}n\geq\frac{d}{2}n^{1-\lambda
/2}>n^{1-\lambda},
\]
completing the proof.
\end{proof}

\section{\label{Degs}Degenerate and splittable graphs}

Proposition \ref{p3} follows from the corollary to the following lemma.

\begin{lemma}
\label{lts}Let $k\geq1,n\geq2$ be integers. For any tree $T_{n}$ of order $n,$
there exists a set $S_{k}\subset V\left(  T_{n}\right)  $ such that
$\left\vert S_{k}\right\vert \leq2^{k+2}-6$ and $\psi\left(  T_{n}%
-S_{k}\right)  \leq2^{-k}n$.
\end{lemma}

\begin{proof}
We shall use induction on $k.$ According to a result from \cite{EFRS88a},
either $\psi\left(  T_{n}-uv\right)  \leq2n/3$ for some $uv\in E\left(
T_{n}\right)  ,$ or $\psi\left(  T_{n}-u\right)  \leq n/3$ for some $u\in
V\left(  T_{n}\right)  .$ Therefore, $\psi\left(  T_{n}-u-v\right)  \leq n/2$
for some vertices $u,v\in V\left(  T_{n}\right)  ,$ implying the lemma for
$k=1$ with $S_{1}=\left\{  u,v\right\}  .$ Assume the lemma holds for $k-1$
and let $S_{k-1}$ be a set such that $\psi\left(  T_{n}-S_{k-1}\right)
\leq2^{-k+1}n.$ For each component $C$ of $T_{n}-S_{k-1}$ with $\left\vert
C\right\vert >2^{-k}n,$ select two vertices $u,v\in V\left(  C\right)  $ such
$\psi\left(  C-u-v\right)  \leq\left\vert C\right\vert /2\leq2^{-k}n.$ Since
there are fewer than $2^{k}$ components $C$ satisfying $\left\vert
C\right\vert >2^{-k}n,$ we deduce that $\left\vert S_{k}\right\vert
<\left\vert S_{k-1}\right\vert +2^{k+1},$ completing the induction step and
the proof.
\end{proof}

\begin{corollary}
Suppose $0<\gamma<1$ is fixed. For every $0<\eta<1,$ every sufficiently large
tree is $\left(  \gamma,\eta\right)  $-splittable.
\end{corollary}

\begin{proof}
Set $k=\left\lceil \log_{2}1/\varepsilon\right\rceil .$ Lemma \ref{lts}
implies that there exists $S\subset V\left(  T_{n}\right)  $ such that
$\left\vert S\right\vert <2^{k+2}-6$ and $\psi\left(  T_{n}-S\right)
\leq2^{-k}n\leq\eta n$. We deduce that $\left\vert S\right\vert <2^{k+2}%
-6<2^{k+2}<8\eta^{-1}<n^{1-\gamma}$ for $n$ large.
\end{proof}

Next we sketch the proofs of Proposition \ref{p4} and \ref{p1}.

\begin{proof}
[\textbf{Proof of Propostion \ref{p4}}]If $\Delta\left(  G\right)  \leq q$
then $\Delta\left(  G^{k}\right)  \leq q^{k};$ hence $\mathcal{F}^{k}$ is
degenerate. Let $\mathcal{F}$ be $\gamma$-crumbling, $G\in\mathcal{F}$ is a
graph of order $n$ and $M\subset V\left(  G\right)  $ is a set such that
$\left\vert M\right\vert <n^{1-\gamma}$ and $\psi\left(  G-M\right)
<\varepsilon n.$ Set%
\[
\left\{  M^{\prime}=v:v\in V\left(  G\right)  ,\text{ there exists }u\in
M\text{ with }dist\left(  u,v\right)  \leq k\right\}  .
\]
If $A$ and $B$ are components of $G-M,$ then $dist\left(  A-M^{\prime
},B-M^{\prime}\right)  \geq2k.$ Therefore, $\psi\left(  G^{k}-M^{\prime
}\right)  <\varepsilon n,$ implying that $\mathcal{F}^{k}$ is $\left(
\gamma/2\right)  $-crumbling.
\end{proof}

\bigskip

\begin{proof}
[\textbf{Proof of Propostion \ref{p1}}]Burr and Erd\H{o}s (\cite{BuEr83},
Lemma 5.4) proved that for every graph $G$ there exists $k\geq1$ such that
every graph of order $n$ homeomorphic to $G$ can be embedded in $P_{n}^{k}$.
This completes the proof, in view of Propositions \ref{p2} and \ref{p4}.
\end{proof}

\bigskip

\begin{proof}
[\textbf{Proof of Propostion \ref{p5}}]Observe that if $G_{1}$ is $q_{1}%
$-degenerate and $G_{2}$ is $q_{2}$-degenerate then $G_{1}\times G_{2}$ is
$\left(  q_{1}+q_{2}\right)  $-degenerate. Also let $G_{1}=G\left(  n\right)
$ be a $\left(  \gamma_{1},\eta_{1}\right)  $-splittable graph and
$G_{2}=G\left(  m\right)  $ be a $\left(  \gamma_{2},\eta_{2}\right)
$-splittable graph. Suppose $m\leq n,$ select $M\subset V\left(  G_{1}\right)
$ with $\left\vert M\right\vert <n^{1-\gamma_{2}}$ such that $\psi\left(
G_{1}-M\right)  <\eta_{1}n.$ Then
\[
\left\vert M\times V\left(  G_{2}\right)  \right\vert =n^{1-\gamma_{1}}%
m\leq\left(  mn\right)  ^{1-\gamma_{1}/2}%
\]
and $\psi\left(  G_{1}\times G_{2}-M\times V\left(  G_{2}\right)  \right)
<\eta_{1}nm.$ Therefore, the graph $G_{1}\times G_{2}$ is $\left(  \gamma
,\eta\right)  $-splittable with $\gamma=\min\left\{  \gamma_{1}/2,\gamma
_{2}/2\right\}  $ and $\eta=\max\left\{  \eta_{1},\eta_{2}\right\}  $
\end{proof}

\bigskip

\begin{proof}
[\textbf{Proof of Propostion \ref{p6}}]Let $\mathcal{F}$ be a $\gamma
$-crumbling family. Suppose $G\in\mathcal{F}$ is a graph of order $n$ and
$M\subset V\left(  G\right)  $ is such that $\left\vert M\right\vert
<n^{1-\gamma}$ and $\psi\left(  G-M\right)  <\eta n.$ Let $\varphi
:G^{\mathbf{k}_{n}}\rightarrow G$ be the homomorphism mapping every vertex to
its ancestor. From the graph $G^{\mathbf{k}_{n}}$ remove the set $M^{\prime
}=\varphi^{-1}\left(  M\right)  .$ If $C$ is a component of $G-M,$ then
$\varphi^{-1}\left(  C\right)  $ is a component of $G^{\mathbf{k}_{n}%
}-M^{\prime}$ and so
\[
\psi\left(  G^{\mathbf{k}_{n}}-M^{\prime}\right)  \leq K\psi\left(
G-M\right)  <K\eta n.
\]
Also,
\[
\left\vert M^{\prime}\right\vert \leq K\left\vert M\right\vert <Kn^{1-\gamma
}<\left(  Kn\right)  ^{1-\gamma/2}%
\]
for $n$ large. Hence, $\left\{  G^{\mathbf{k}_{n}}\right\}  $ is a $\left(
\gamma/2\right)  $-crumbling family.
\end{proof}

\section{\label{Discon}Disproof of Conjecture \ref{con2}}

In this section we shall prove the following result.

\begin{theorem}
For $n$ sufficiently large, almost all connected $100$-regular graphs of order
$n$ are not $3$-good.
\end{theorem}

Our idea is a refinement of the main idea in \cite{Bra96}; however to simplify
the presentation, we use newer, more powerful results.

Define a $2$-coloring $E\left(  K_{2n-1}\right)  =E\left(  R\right)  \cup
E\left(  B\right)  $ as follows. Partition $V(K_{2n-1})=[2n-1]$ into five sets
$V_{1},\ldots,V_{5}$ so that $\left\vert V_{1}\right\vert \leq\ldots
\leq\left\vert V_{5}\right\vert \leq\left\vert V_{1}\right\vert +1;$ thus,
each set has $\left\lfloor (2n-1)/5\right\rfloor $ or $\left\lceil
(2n-1)/5\right\rceil $ vertices. Set $E\left(  R\right)  =\left\{  uv:u\in
V_{i},\text{ }v\in V_{j},\text{ }i-j\equiv\pm1\pmod{5}\right\}  $ and let all
other edges belong to $E\left(  B\right)  $. Clearly, the graph $R$ is $K_{3}%
$-free. We claim that, for $n$ sufficiently large, $G\nsubseteq B$ for almost
all connected $100$-regular graphs $G$ of order $n$. To prove this claim we
need first a proposition.

\begin{proposition}
Every subgraph of $B$ of order $n$ contains two disjoint sets sets $X$ and $Y$
with $\left\vert X\right\vert \left\vert Y\right\vert \geq n^{2}/25-O\left(
n\right)  $ and $e_{B}\left(  X,Y\right)  =0.$
\end{proposition}

\begin{proof}
Let $q(n)$ be the largest integer such that every $n$-element subset of
$V(K_{n-1})=[2n-1]$ induces a complete bipartite subgraph of size $q(n)$ in
$R$. We shall prove that%
\[
q(n)>\frac{n^{2}}{25}-O(n),
\]
implying the desired result.

Let $X$ be an $n$-element subset of $[2n-1]$, and set $X_{i}=X\cap V_{i}$ for
$1\leq i\leq5$. We may assume that $\left\vert X_{5}\right\vert =\max
\limits_{i}\left\vert X_{i}\right\vert $. Note that $X$ induces two complete
bipartite graph in $R$ - one with parts $X_{5}$ and $X_{1}\cup X_{2}$ and
another one with parts $X_{2}$ and $X_{3}.$ Since $\sum_{i}\left\vert
X_{i}\right\vert =n$, either $\left\vert X_{1}\right\vert +\left\vert
X_{4}\right\vert +\left\vert X_{5}\right\vert \geq n/2$ or $\left\vert
X_{2}\right\vert +\left\vert X_{3}\right\vert \geq n/2$. We consider each of
these two possibilities in turn. If $\left\vert X_{1}\right\vert +\left\vert
X_{4}\right\vert +\left\vert X_{5}\right\vert \geq n/2,$ then $3\left\vert
X_{5}\right\vert \geq n/2$ and $\left\vert X_{5}\right\vert \leq\left\lceil
\left(  2n-1\right)  /5\right\rceil $. Since $x(n/2-x)$ is a concave function
of $x,$ its minimum over $[a,b]$ is $\min\{a(n/2-a),b(n/2-b)\}$. Thus, the
size of the complete bipartite graph with parts $X_{5}$ and $X_{1}\cup X_{2}$
is at least
\[
\left\vert X_{5}\right\vert \left(  n/2-\left\vert X_{5}\right\vert \right)
\geq\min\left\{  \frac{n}{6}\left(  \frac{n}{2}-\frac{n}{6}\right)
,\left\lceil \frac{2n-1}{5}\right\rceil \left(  \frac{n}{2}-\left\lceil
\frac{2n-1}{5}\right\rceil \right)  \right\}  =\frac{n^{2}}{25}-O(n).
\]
Suppose $\left\vert X_{2}\right\vert +\left\vert X_{3}\right\vert \geq n/2$
and assume that $\left\vert X_{2}\right\vert \geq\left\vert X_{3}\right\vert
$. Then $n/4\leq\left\vert X_{2}\right\vert \leq\left\lceil \left(
2n-1\right)  /5\right\rceil $. As before we find that the size of the complete
bipartite subgraph with parts $X_{2}$ and $X_{3}$ is at least $n^{2}/25-O(n)$,
completing the proof.
\end{proof}

Recently Friedman \cite{Fri04} confirmed a conjecture of Alon, proving the
following result.

\begin{fact}
\label{Fried}For even $d\geq4$ and every $\varepsilon>0$, the second singular
value $\sigma_{2}$ of almost all $d$-regular graphs satisfies%
\[
\sigma_{2}\leq2\sqrt{d-1}+\varepsilon.
\]

\end{fact}

Earlier, Robinson and Wormald \cite{RoWo94} proved that for $d\geq3,$ almost
all $d$-regular graphs are Hamiltonian. Therefore, we have the following
simple corollary.

\begin{fact}
\label{RoWo}For $d\geq3,$ almost every $d$-regular graph is connected.
\end{fact}

We need also the following statement, generally known as the \textquotedblleft
Expander mixing lemma\textquotedblright, (for a proof see \cite{KrSu06}, p. 11).

\begin{fact}
\label{EML}For every $d$-regular graph $G$ of order $n$ and every nonempty
sets $X,Y\subset V\left(  G\right)  ,$%
\[
\left\vert e\left(  X,Y\right)  -\frac{d}{n}\left\vert X\right\vert \left\vert
Y\right\vert \right\vert \leq\sigma_{2}\left(  G\right)  \sqrt{\left\vert
X\right\vert \left\vert Y\right\vert }.
\]

\end{fact}

Facts\ \ref{Fried} and \ref{RoWo} imply that almost every $100$-regular graph
$G$ is connected and satisfies%
\[
\sigma_{2}\left(  G\right)  \leq2\sqrt{d-1}+\varepsilon.
\]
If such a graph of sufficiently large order is $3$-good, then Proposition
implies that $G$ contains two disjoint sets $X$ and $Y$ such that $\left\vert
X\right\vert \left\vert Y\right\vert \geq n^{2}/25-O\left(  n\right)  $ and
$e\left(  X,Y\right)  =0.$ Hence,%
\[
\frac{100}{n}\left\vert X\right\vert \left\vert Y\right\vert \leq\sigma
_{2}\left(  G\right)  \sqrt{\left\vert X\right\vert \left\vert Y\right\vert }%
\]
and so,%
\[
20\left(  1+o\left(  1\right)  \right)  \leq\frac{100}{n}\sqrt{\left\vert
X\right\vert \left\vert Y\right\vert }\leq\sigma_{2}\left(  G\right)
\leq2\sqrt{99}+\varepsilon,
\]
a contradiction for large $n$ and $\varepsilon$ sufficiently small.

\end{document}